\documentclass{amsart}
\usepackage{amsmath}
\usepackage{amsfonts}
\usepackage{eufrak}
\usepackage{amsthm}
\usepackage{amssymb}
\usepackage[all]{xy}
 \SelectTips{cm}{}
\title{Uniqueness of $A_\infty$-structures and Hochschild cohomology}
\author{Constanze Roitzheim}
\author{Sarah Whitehouse}
\email{Constanze.Roitzheim@glasgow.ac.uk, S.Whitehouse@sheffield.ac.uk}
\address{C. Roitzheim \\ University of Glasgow \\ Department of Mathematics\\ University Gardens\\ Glasgow G12 8QW,
UK}
\address{S. Whitehouse\\ School of Mathematics and Statistics\\
University of Sheffield \\ Hicks Building\\ Hounsfield Road\\
Sheffield S3 7RH, UK}
\thanks{Supported by EPSRC grant EP/E022618/1}
\date{$8^\text{th}$ October 2010}

\DeclareMathOperator{\Hom}{Hom} 
\DeclareMathOperator{\hh}{HH} 

\theoremstyle{definition}
\newtheorem{ddd}{Definition}[section]
\newtheorem{ppp}[ddd]{Proposition}
\newtheorem{lem}[ddd]{Lemma}
\newtheorem{theorem}[ddd]{Theorem}
\newtheorem{ccc}[ddd]{Corollary}
\newtheorem*{rrr}{Remark}
\newtheorem*{nnn}{Notation}

\newtheorem*{ex}{Example}

\newcommand{\bc}[1]{[\![#1]\!]}

\begin{document}

\begin{abstract}
Working over a commutative ground ring, we establish a Hochschild cohomology
criterion for uniqueness of derived $A_\infty$-algebra structures in the sense of Sagave.
 We deduce a Hochschild cohomology criterion
for intrinsic formality of a differential graded algebra. This
generalizes a classical result of Kadeishvili for the case of a graded algebra over
a field. 
\end{abstract}

\maketitle

\section*{Introduction}
\label{sec:intro}

$A_\infty$-structures were introduced by Stasheff~\cite{Sta63} in
the early 1960s in the study of topological spaces with products.
They are now known to arise widely in algebra, geometry and
mathematical physics, as well as topology.

We are interested in questions of formality and intrinsic formality
for differential graded algebras. Thus we would like to establish
conditions under which two differential graded algebras with the
same homology are quasi-isomorphic. This has been studied by Keller
and others in the case where the ground ring $k$ is a field. It is
related to the existence of different $A_\infty$-structures on a
minimal model of the differential graded algebra.

An important structural result of Kadeishvili~\cite{Kad79} proves 
the existence of minimal models of differential graded algebras 
over a field while another classical theorem by Kadeishvili~\cite{Kad88} 
gives a criterion for uniqueness of certain minimal models using Hochschild cohomology.

\medskip
For the applications we have in mind, which are related to
rigidity of the model category structures arising in stable homotopy
theory, we will be interested in working over local rings rather
than fields. When working with a commutative ground ring rather than
a field, one has to work with derived $A_\infty$-algebras as in the
world of ``classical'' $A_\infty$-algebras, a differential graded
algebra might not have a minimal model if its homology is not
projective. The theory of derived $A_\infty$-algebras was developed
by Sagave in~\cite{Sag10}. He describes the notion of a minimal
model for a differential graded algebra $A$ over a commutative
ground ring by giving a projective resolution of the homology of $A$
that is compatible with the existing $A_\infty$-structure on $A$.

\medskip
Our main result is Theorem~\ref{derivedHHcriterion} which extends Kadeishvili's uniqueness 
theorem to derived $A_\infty$-algebras. For this we develop a new notion of Hochschild cohomology. 
After some further work we again
obtain a Hochschild cohomology criterion for intrinsic formality of
a differential graded algebra over a commutative ring rather than a
field, Theorem~\ref{th:derivedformality}.

In the subsequent sections we return to classical $A_\infty$-algebras and derive
some further generalizations of Kadeishvili's uniqueness criterion. 
The first of these is Theorem~\ref{th:HHcriterion} which studies uniqueness of 
an $A_\infty$-structure on a fixed differential graded algebra. The other, Theorem~\ref{uniquenessmassey}, 
discusses differential graded algebras with fixed Massey products on their homology. 

\medskip
An alternative approach is developed by Dugger and Shipley.
In~\cite[Section 3]{DugShi07} they consider the classification of
quasi-isomorphism types of differential graded algebras with given
homology. They do this by building differential graded algebras up
degreewise via a theory of Postnikov sections and $k$-invariants. To
do so requires working with bounded below differential graded
algebras, a restriction which does not apply to our methods. The
$k$-invariants live in derived Hochschild cohomology groups of the
Postnikov sections with coefficients in the next homology group of
the differential graded algebra being built. Their work does not
consider $A_\infty$-structures and although also formulated in terms
of Hochschild cohomology, there does not seem to be a very direct
relationship between their methods and ours. However, we are 
going to put some of their examples in context throughout our paper.
\medskip

This paper is organized as follows. In Section~\ref{sec:Ainf} we
recall basic definitions relating to $A_\infty$-algebras and
Hochschild cohomology. In Section
\ref{sec:dAinf} we recall Sagave's construction of derived
$A_\infty$-algebras and his results about minimal models. 
This section also introduces the Lie algebra structure which leads to 
the definition of Hochschild cohomology of a certain class of derived 
$A_\infty$-algebras in Section~\ref{sec:deruniq}. At the end of Section~\ref{sec:deruniq} 
we show that the vanishing of certain Hochschild cohomology groups
gives a sufficient condition for the existence of a unique derived 
$A_\infty$-structure on a fixed underlying object. In Section~\ref{sec:e2inv}
we deduce the criterion for intrinsic formality of 
differential graded algebras over a commutative ground ring.
Finally, 
in Sections \ref{sec:unique} and \ref{sec:massey} we discuss the 
previously mentioned analogues of these results
 for 
classical $A_\infty$-structures. A short appendix is devoted to sign issues.

\bigskip

{\bf Acknowledgments.} We would like to thank Andy Baker, David
Barnes and Fernando Muro for motivating comments and suggestions.
Further thanks go to Steffen Sagave for patiently answering
questions about sign conventions.

\section{A quick review of $A_\infty$-algebras}\label{sec:Ainf}

We assume that the reader is familiar with the basic 
definitions regarding $A_\infty$-algebras and Hochschild cohomology, 
but we are going to recall some of them in this section to establish 
notation and assumptions. We are going to be very brief with this; 
the explicit formulas and definitions regarding derived 
$A_\infty$-algebras given in the later Sections~\ref{sec:dAinf} 
and~\ref{sec:deruniq} specialize to the case of ``classical'' 
$A_\infty$-algebras. For greater detail we refer to Keller's introductory paper~\cite{Kel01}.

The notion of an $A_\infty$-algebra arose with the study of loop
spaces in topology and has since become an increasingly important
and powerful subject in algebraic topology and homological algebra. Roughly speaking,
$A_\infty$-algebras are not necessarily associative algebras with
given maps for ``multiplying'' $n$ elements for each $n$, unlike in
the case  of associative algebras where one knows how to multiply
$n$ elements from knowing how to multiply two elements.

\subsection{Basic definitions}

\medskip In Sections~\ref{sec:Ainf} and \ref{sec:massey} of this
paper, $k$ will denote a field of characteristic not equal to $2$.
In Sections~\ref{sec:dAinf} to~\ref{sec:unique} we will allow $k$ to
be a commutative ring rather than a field. Note that in fact
Sections~\ref{sec:Ainf} and \ref{sec:massey} do not require a ground
field as long as all $k$-modules in question are projective.

All unadorned tensor products are over $k$. All graded objects will
be $\mathbb{Z}$-graded unless stated otherwise. Our convention for the degree of a map $f$
is as follows: a map of graded $k$-vector spaces $f:A\to B$ of
degree $i$ consists of a sequence of maps $f^n:A^n\to B^{n+i}$.
(Later this will be called the internal degree and there will also
be a notion of cohomological or external degree.) We often
abbreviate `differential graded algebra' to dga.

\begin{ddd}
Let $A=\bigoplus\limits_{n \in \mathbb{Z}} A^n$ be a graded
$k$-vector space. An {\it $A_\infty$-structure} on $A$ is a sequence
of $k$-linear maps
\[
m_j: A^{\otimes j} \longrightarrow A \quad \text{for $j\geq 1$}
\]
of degree $2-j$ satisfying the equation
\[
\sum_{n=r+s+t} (-1)^{rs+t} m_{1+r+t}(1^{\otimes r} \otimes m_s
\otimes 1^{\otimes t})=0
\]
for each $n\ge 1$. An {\it $A_\infty$-algebra} is a graded
$k$-vector space $A$ together with an $A_\infty$-structure on $A$.

\medskip Further all $A_\infty$-algebras are assumed to be strictly
unital; c.f.~Definition~\ref{derivedainfstructure}.
We are using the
sign convention of Sagave~\cite[(2.6)]{Sag10} and of
Lef\`evre-Hasegawa~\cite[1.2.1.2]{LH03} rather than of
Keller~\cite{Kel01}.
\end{ddd}

Note that we are applying the Koszul sign rule when applying such
formulas to elements:
\[
(f \otimes g)(x \otimes y)=(-1)^{|g| |x| }f(x) \otimes g(y).
\]

In particular, this definition gives us
\[
m_1 m_1 = 0,
\]
i.e.~$m_1$ is a differential on $A$. It also yields the following special
cases: if $m_k=0$ for all $k \neq 2$, then $A$ is simply a graded
associative algebra. If $m_k=0$ for $k \ge 3$, then $A$ is a
differential graded algebra.

\medskip
There are also notions of morphism and quasi-isomorphism of $A_\infty$-algebras;
these are special cases of Definitions~\ref{derivedainfmorphism} and~\ref{def:e2equivalence}.

\begin{nnn}
We sometimes write an $A_\infty$-structure as a formal infinite sum,
i.e.\[ m= m_1 + m_2 + \cdots .\] Note that all infinite sums in this
paper are finite in every degree.
\end{nnn}

\subsection{Hochschild cohomology and Lie structure}\label{liealgebrastructure}

Hochschild cohomology is a very powerful tool in many areas around
algebra and topology, from relations to the geometry of loop spaces
to deformation theory of algebras and realizability questions in
topology. The definition of Hochschild cohomology of associative graded 
algebras can be extended to a definition of Hochschild cohomology 
of $A_\infty$-algebras. A convenient way of doing this is using a Lie algebra structure on the bigraded $k$-vector space 
\[
C^{n,m}(A,A)=\Hom_k^m(A^{\otimes
n},A)=\prod\limits_i\Hom_k((A^{\otimes n})^{i},A^{i+m}),
\]
where $n\in\mathbb{N}, m\in\mathbb{Z}$ and $A$ is a graded $k$-vector space.

Explicitly, for $f\in C^{n,k}(A,A)$ and $g \in C^{m,l}(A,A)$ the Lie bracket is given by
\begin{align}
[f,g] &= \sum\limits_{i=0}^{n-1} (-1)^{(n-1)(m-1)+(n-1)l+i(m-1)}
f(1^{\otimes i} \otimes g \otimes 1^{\otimes n-i-1}) \nonumber\\
 & -(-1)^{(n+k-1)(m+l-1)}\sum\limits_{i=0}^{m-1}
(-1)^{(m-1)(n-1)+(m-1)k+i(n-1)} g(1^{\otimes i} \otimes f \otimes
1^{\otimes m-i-1})  \nonumber
\end{align}
which lies in $C^{n+m-1,l+k}(A,A)$.
This gives $C^{*,*}(A,A)$ the structure of a graded Lie algebra,
where the grading is by total degree shifted by $1$; see
e.g. \cite[Section 2]{FiaPen02}, \cite[Section 1]{Get94},
\cite{Ger63} or \cite{PenSch95}.
Note that the formula given in some
of the references has signs arising from the Koszul rule because it
is given evaluated on elements rather than as a formula of
morphisms. For details on how this formula arises, see Section \ref{trigraded} and the Appendix.

\begin{lem}\label{zerobracket}
Let $m \in C^{*,*}(A,A)$ of total degree 2. Then $m$ is an
$A_\infty$-structure on $A$ if and only if $[m,m]=0$. Further, for such
$m$,
\[
D:= [m,-]: C^{*,*}(A,A) \longrightarrow C^{*,*}(A,A)
\]
is a differential on $C^{*,*}(A,A)$, i.e. $D$ raises total
degree by 1 and satisfies $D \circ D = 0$.
\end{lem}

\begin{proof}
The first claim follows immediately from the bracket formula and the fact that
$2$ is invertible. The
fact that $D\circ D=0$ is an immediate consequence of the graded
Jacobi identity, while the total degree of $D$ can be computed
directly.
\end{proof}

\begin{ddd}\label{def:HHclassical}
Let $A$ be an $A_\infty$-algebra with $A_\infty$-structure $m$. Then
the {\it Hochschild cohomology of the $A_\infty$-algebra A} is
defined as
\[
\hh^{*}(A,A)=H^{*}\left(\prod\limits_{i} C^{i,*-i}(A,A), [m,-]\right).
\]
\end{ddd}

For this, see, for example,~\cite[\S 5]{PenSch95}. If $A$ is an
associative algebra (i.e. $m=m_2$), a direct computation using the above
definitions shows this recovers the usual definition of the
Hochschild cohomology of associative algebras, i.e. for $f\in C^{n,k}(A,A)$,
\[
[m_2,f]=(-1)^k \Big(m_2(1\otimes f)
+\sum\limits^{n-1}_{i=0}(-1)^{i+1}f(1^{\otimes i} \otimes m_2
\otimes 1^{\otimes n-1-i})  +(-1)^{n+1} m_2(f \otimes 1) \Big).
\]

The grading in Definition~\ref{def:HHclassical} refers to the total degree. 
In the case of an associative algebra the differential 
\[
[m_2,-]: C^{*,*}(A,A) \longrightarrow C^{*+1,*}(A,A)
\]
preserves internal degree so we can split the total degree of the Hochschild 
cohomology into the cohomological degree and the internal degree. We 
denote the bigraded Hochschild cohomology in this special case by 
$\hh^{*,*}_{alg}(A,A)$. 
For a general $A_\infty$-algebra, we do not have a bigrading, but we can introduce a filtration,
see Definition~\ref{def:filtration}.

\medskip
For $A$ a dga, the definition can be interpreted in terms of
bicomplexes. The dga $A$ has differential $m_1$ and multiplication
$m_2$.  The bigraded module $C^{*,*}(A,A)$ becomes a bicomplex by
taking
\[
d^v:=[m_1,-]: C^{*,*} \longrightarrow C^{*,*+1}
\]
to be the vertical differential and
\[
d^h:=[m_2,-]: C^{*,*} \longrightarrow C^{*+1,*}
\]
to be the horizontal differential. 
The condition
\[
[m_1+ m_2, m_1 + m_2] =0
\]
translates into $(d^v)^2=0$, $(d^h)^2=0$ and $d^v d^h + d^h d^v =
0$, which are exactly the conditions for $C^{*,*}(A,A)$ to be a
bicomplex \cite[1.2.4]{Wei94}.

\subsection{Minimal models and uniqueness}

We now recall a definition and theorem about minimal models of
$A_\infty$-algebras. It relates differential graded algebras to $A_\infty$-structures on their homology.

\begin{ddd}
An $A_\infty$-algebra is called {\it minimal} if $m_1 = 0$.
\end{ddd}

Over a field, one can replace any $A_\infty$-algebra by a quasi-isomorphic
minimal one which gives a very convenient way to describe a
quasi-isomorphism class of an $A_\infty$-algebra. We are 
particularly interested in the special case of differential graded
algebras.

\pagebreak

\begin{theorem}[Kadeishvili]\label{Kadeishviliminimal}
Let $A$ be a differential graded algebra over a field $k$, and let $H^*(A)$ be its homology
module. Then $H^*(A)$ has an $A_\infty$-structure such that
\begin{itemize}
\item $m_1=0$ and the multiplication $m_2$ is induced by the multiplication on $A$,
\item there is a morphism of $A_\infty$-algebras $f: H^*(A) \longrightarrow A$ such that $f_1$ is a quasi-isomorphism.
\end{itemize}
This $A_\infty$-algebra $H^*(A)$ is called the {\it minimal model of
$A$}.
\end{theorem}

For more details, see~\cite{Kad79}. Note that the theorem states
in particular that the minimal model $H^*(A)$ is quasi-isomorphic to
$A$ as an $A_\infty$-algebra.
\smallskip
This is useful in combination with a uniqueness result in \cite{Kad88}

\begin{ddd}
We say that an $A_\infty$-structure $m$ is \emph{trivial} if $m_n=0$
for $n \ge 3$.
\end{ddd}

\begin{theorem}[Kadeishvili]\label{Kadeishviliuniqueness}
Let $C$ be a graded $k$-algebra with multiplication $\mu$. If 
\[
\hh^{n,2-n}_{alg}(C,C)=0 \quad\quad\mbox{for $n\ge 3$},
\]
then every $A_\infty$-structure on $C$ with $m_1=0$ and $m_2=\mu$ is quasi-isomorphic to the trivial one.
\end{theorem}

We can reformulate this in terms of formality of dgas. We recall the following
standard definitions.

\begin{ddd}
\label{def:formal}

\begin{enumerate}
\item A dga $A$ is \emph{formal} if it is quasi-isomorphic to its homology $H^*(A)$ regarded
as a dga with trivial differential.

\item A dga $A$ is \emph{intrinsically formal} if any other dga $A'$ such that $H^*(A)\cong H^*(A')$
as associative algebras is quasi-isomorphic to $A$.
\end{enumerate}
\end{ddd}

If a dga is intrinsically formal then it is formal, but the converse need not hold.
For example, in~\cite[Example 3.15]{DugShi07}, it is shown that there are two quasi-isomorphism
types of dgas with homology an exterior algebra over $\mathbb{F}_p$ on an even degree generator.
The trivial one is therefore formal but not intrinsically formal.
\smallskip

Using Theorem \ref{Kadeishviliuniqueness} for the case $C=H^*(A)$ yields the following.

\begin{ccc}\label{uniquedga}
Let $A$ be a dga and $H^*(A)$ its homology algebra. Suppose that
\[
\hh^{n,2-n}_{alg}(H^*(A),H^*(A))=0 \quad\quad\mbox{for } n \ge 3.
\]
Then $A$ is intrinsically formal.
\end{ccc}

In Section~\ref{sec:unique}, we will recover these results as special cases of our derived versions.

\section{Derived $A_\infty$-algebras}\label{sec:dAinf}

To work with Kadeishvili's minimal models and to establish the
uniqueness theorems, one has to assume all dgas as well as their
homology algebras to be degreewise projective, hence the assumption of a ground field. However, there are
important examples arising from homotopy theory where projectivity
cannot be guaranteed. In 2008, Sagave introduced the notion of
derived $A_\infty$-algebras, providing a framework for not
necessarily projective modules over an arbitrary commutative ground
ring~\cite{Sag10}.

First of all, we recall some definitions and results about derived
$A_\infty$-algebras; we refer to Sagave's paper for the finer
technical details.

\medskip
The basic idea is to introduce degreewise projective resolutions for
an $A_\infty$-algebra that are compatible with the
$A_\infty$-structure. This will introduce another internal grading.

\subsection{Definitions, conventions and known
results}

All definitions and results in this subsection have been developed
by Sagave in~\cite{Sag10} and we refer to his paper for technical
details.

Let $k$ be a commutative ring and let
$A$ be
an {\it $(\mathbb{N}$,$\mathbb{Z})$-bigraded} $k$-module, i.e.~$A =
\bigoplus\limits_{i\in \mathbb{N},j\in\mathbb{Z}} A^{j}_i$. A
morphism of bigraded $k$-modules $f: A \longrightarrow B$ of
bidegree $(s,t)$ is a sequence of maps of $k$-modules $f: A^{j}_i
\longrightarrow B_{i-s}^{j+t}$ for all $i \in \mathbb{N}$ and $j \in
\mathbb{Z}$. Again, we follow the Koszul sign convention: for $g$ a
morphism of bidegree $(s,t)$ and $x$ an element of bidegree $(i,j)$,
we have
$$(f \otimes g)(x \otimes y) = (-1)^{is + jt}f(x) \otimes g(y).$$
The homological (subscript) bidegree is called the {\it horizontal
bidegree} and the cohomological (superscript) bidegree is called the
{\it vertical bidegree}.

Throughout the rest of the paper we also assume that all bigraded modules have no 2-torsion.

\begin{ddd}\cite[Definition 2.1]{Sag10}\label{derivedainfstructure}
A {\it derived $A_\infty$-structure} (or {\it $dA_\infty$-structure}
for short) on an $(\mathbb{N}$,$\mathbb{Z})$-bigraded $k$-module $A$
consists of $k$-linear maps
\[
m_{ij}: A^{\otimes j} \longrightarrow A
\]
of bidegree $(i,2-(i+j))$ for each $j\ge 1$, $i \ge 0$, satisfying the equation
\begin{equation}\label{dobjectequation}
\sum\limits_{\substack{ u=i+p, v=j+q-1 \\ j=1+r+t}} (-1)^{rq+t+pj}
m_{ij} (1^{\otimes r} \otimes m_{pq} \otimes 1^{\otimes t}) = 0
\end{equation}
for all $u\ge 0$ and $v \ge 1$. A {\it $dA_\infty$-algebra} is a
bigraded $k$-module together with a $dA_\infty$-structure.

A $dA_\infty$-algebra $A$ is called {\it strictly unital} if there
is a unit map $\eta: k \longrightarrow A$ such that
\begin{itemize}
\item $m_{01}(\eta)=0$,
\item $m_{02}(\eta \otimes 1) = 1 = m_{02}(1 \otimes \eta)$,
\item $m_{ij}(1^{\otimes r-1}\otimes \eta\otimes 1^{\otimes j-r})=0$ for $i+j \ge
3$, $1 \le r \le j$.
\end{itemize}
\end{ddd}

\smallskip

From now on, all $dA_\infty$-algebras are assumed to be strictly
unital.
\begin{rrr}A $dA_\infty$-algebra concentrated in horizontal degree 0
(and hence with $m_{ij}=0$ for all $i \neq 0$) is the same as an
$A_\infty$-algebra.

A $dA_\infty$-algebra with $m_{ij}=0$ except $m_{01}$ and $m_{11}$
is just a bicomplex (with a different sign convention to that
encountered earlier) with horizontal differential $m_{11}$ and
vertical differential $m_{01}$ as the definition in this case forces
$m_{11}m_{11}=0$, $m_{01}m_{01}=0$ and
$m_{01}m_{11}-m_{11}m_{01}=0$.
\end{rrr}

\begin{ddd}
A {\it bidga} is a monoid in the category of bicomplexes; equivalently, a bidga is a
$dA_\infty$-algebra with $m_{ij}=0$ for $i+j \ge 3$.
(See~\cite[Definition 2.10 and Remark 2.11]{Sag10}.)
\end{ddd}

\begin{ddd}\cite[Definition 2.5]{Sag10}\label{derivedainfmorphism}
Let $A$ and $B$ be $dA_\infty$-algebras with $dA_\infty$-structures
$m$ and $\overline{m}$, respectively. A {\it morphism of
$dA_\infty$-algebras} $f: A \longrightarrow B$ consists of a family
of $k$-module maps
\[
f_{st}: A^{\otimes t} \longrightarrow B
\]
of bidegree $(s,1-(s+t))$ satisfying
\begin{align}
\sum\limits_{\substack{u=i+p, v=j+q-1 \\ j=1+r+t }} (-1)^{rq+t+pj}
f_{ij} (1^{\otimes r} \otimes m_{pq} \otimes 1^{\otimes t}) =
\sum\limits_{\substack{ u= i + p_1 + \cdots + p_j \\ v= q_1+\cdots+ q_j}}
(-1)^{\epsilon} \quad \overline{m}_{ij} (f_{p_1 q_1} \otimes \cdots
\otimes f_{p_j q_j}) \label{dmorphequation}
\end{align}
for all $u \ge 0$ and $v \ge 1$. Here,
\[
\epsilon = u + \sum\limits_{w=1}^{j-1} \left( jp_w + w(q_{j-w} -
p_w)+q_{j-w}\left(\sum\limits_{s=j-w+1}^{j}p_s + q_s\right)\right).
\]
For strictly unital $dA_\infty$-algebras, morphisms are required to
satisfy the unit conditions $f_{01}\eta=\overline{\eta}$ and
$f_{ij}(1^{\otimes r-1} \otimes \eta\otimes 1^{\otimes j-r})=0$ for
$i+j \ge 2$ and $1 \le r \le j$.
\end{ddd}

Recall that a quasi-isomorphism of $A_\infty$-algebras is a morphism
of $A_\infty$-algebras that induces a quasi-isomorphism of complexes
with respect to $m_1$. In the case of $dA_\infty$-algebras, the role
of the quasi-isomorphisms is played by the so-called
$E_2$-equivalences. These are the morphisms that induce an
isomorphism of $E_2$-terms of the spectral sequence computing the
homology of the total complex of a bicomplex, see
\cite[2.12]{Mcc01}.

\begin{nnn}
The equations defining a $dA_\infty$-structure include
$m_{01}m_{01}=0$. For a $dA_\infty$-algebra $A$ let $H_{ver}^*$
denote its homology with respect to the {\it vertical differential}
$m_{01}$. The map $m_{01}$ is called the vertical differential
because it raises the vertical degree.
\end{nnn}

Since the equations defining a $dA_\infty$-structure also include
$m_{21}m_{01}-m_{11}m_{11}+m_{01}m_{21}=0$, it follows that the map
$m_{11}$ becomes a differential in horizontal direction on the
bigraded module $H_{ver}^*(A)$, so we can form
$H_{hor}^*(H_{ver}^*(A))=H^*(H_{ver}^*(A),m_{11}).$

\begin{ddd}\label{def:e2equivalence}
A morphism $f: A \longrightarrow B$  of $dA_\infty$-algebras is
called an {\it $E_2$-equivalence} if
$$H_{hor}^*(H_{ver}^*(f_{01}))$$ is an isomorphism of $k$-modules,
c.f. \cite[Definition 2.19]{Sag10}.
\end{ddd}

We would like to extend some applications of $A_\infty$-algebras to
differential graded algebras that are not necessarily projective
over the ground ring $k$ or whose homology is not projective. The
problem we encounter is that not all differential graded algebras
possess a minimal model as an $A_\infty$-algebra. However, Sagave
showed that dgas have reasonable minimal models in the world of
$dA_\infty$-algebras. For this, one has to apply a special
projective resolution.

\begin{ddd}\cite[Definition 3.1]{Sag10}
Let $A$ be a graded algebra. A {\it termwise $k$-projective resolution}
of $A$ is a termwise $k$-projective bidga $P$ with
$m_{01}=0$ together with an $E_2$-equivalence $P \longrightarrow A$.
\end{ddd}

\begin{ddd}\cite[Definition 3.2]{Sag10}
Let $A$ be a dga. A {\it $k$-projective $E_1$-resolution of $A$} is
a bidga $B$ together with an $E_2$-equivalence $B \longrightarrow A$
such that $H_{ver}^{st}(B)$ is projective for each bidegree.
Further, the map $k \longrightarrow H^{00}_{ver}(B)$ induced by the
unit $k \longrightarrow B$ is required to split as a $k$-module map.
\end{ddd}

Thus a $k$-projective $E_1$-resolution of a dga $A$ induces
a termwise $k$-projective resolution of the graded homology algebra 
of $A$.

Sagave then proceeds to show that a $k$-projective $E_1$-resolution
is unique up to $E_2$-equivalence.

\begin{theorem}\cite[Theorem 3.4]{Sag10}
Every dga $A$ over $k$ admits a $k$-projective $E_1$-resolution. Two
such resolutions can be related by a zig-zag of $E_2$-equivalences
between $k$-projective $E_1$-resolutions.
\end{theorem}

\begin{ddd}
A $dA_\infty$-algebra is called {\it minimal} if $m_{01}=0$.
\end{ddd}

\begin{theorem}\cite[Theorem 1.1]{Sag10}\label{derivedminimalmodel}
Let $A$ be a dga over $k$. Then there is a degreewise $k$-projective
$dA_\infty$-algebra $E$ together with an $E_2$-equivalence $E
\longrightarrow A$ such that
\begin{itemize}
\item $E$ is minimal,
\item $E$ is well-defined up to $E_2$-equivalence,
\item together with the differential $m_{11}$ and the multiplication
$m_{02}$, $E$ is a termwise $k$-projective resolution of the graded algebra
$H^*(A)$.
\end{itemize}
\end{theorem}

To prove this, Sagave starts with a $k$-projective $E_1$-resolution
$E \longrightarrow A$. He then shows that the vertical homology
$H^*_{ver}(E)$ admits a $dA_\infty$-structure satisfying the claims
of the theorem.

However, not every termwise projective resolution of $H^*(A)$ admits
such a structure \cite[Remark 4.14.]{Sag10}. For example, consider
the dga over $\mathbb{Z}$
\[
A= \mathbb{Z}[e]\big/(e^4), \quad\quad \partial(e)=p, \quad \quad
|e|=-1,
\]
also examined by Dugger and Shipley in \cite[Example
3.13]{DugShi07}. The bidga
\[
C= \mathbb{Z}\left<a,b\right>\big/(a^2, b^2, ab-ba), \quad\quad
|a|=(1,0), |b|=(0,-2), \quad\quad m_{11}(b)=p
\]
is a termwise projective resolution of
$H^*(A)=\Lambda_{\mathbb{Z}/p}([e^2])$, but there is no
$dA_\infty$-structure on $C$ admitting an $E_2$-equivalence $C
\longrightarrow A$. (For example, equation (\ref{dmorphequation})
for $(u,v)=(2,2)$ forces $m_{22}(b\otimes b)\equiv \pm 1 \bmod p$
whereas equation (\ref{dobjectequation}) for $(u,v)=(2,3)$ forces
$m_{22}(b\otimes b)\equiv 0 \bmod p$.)

\begin{ddd}
Let $A$ and $E$ be as in Theorem~\ref{derivedminimalmodel}. Such an
$E$ is called a {\it minimal model} of $A$.
\end{ddd}

\begin{rrr}
Note that in the context of Theorem~\ref{derivedminimalmodel}, the
underlying $k$-module of the minimal model $E$ together with the
differentials $m_{01}$ and $m_{11}$ and the multiplication $m_{02}$
form a bidga.
\end{rrr}

\subsection{Lie algebra structure on
$C^{*,*}_*(A,A)$}\label{trigraded}

We would like to establish a reasonable notion of Hochschild
cohomology for $dA_\infty$-algebras. In
order to give a simple description, it is our goal to describe the
Hochschild cohomology in terms of a graded Lie algebra structure.

\medskip Let $A$ be a $(\mathbb{N},\mathbb{Z})$-bigraded module without 2-torsion over
a commutative ring. Define
$$C^{n,i}_k(A,A) = \prod\limits_{u,v} \Hom((A^{\otimes n})_{u}^{v},
A_{u-k}^{v+i}).$$

We are going to define a Lie algebra structure on $C^{*,*}_*(A,A)$
generalizing Section \ref{liealgebrastructure}. First of
all, we define a bracket operation that is not a Lie bracket.
Then we are going to introduce a shift operation on elements of
$C^{*,*}_*(A,A)$ and then define the actual Lie bracket using this
shift and the previously defined bracket operation.

For $f \in C^{n,i}_k(A,A)$ and $g \in C^{m,j}_l(A,A)$ we now define

\begin{align}
\bc{f,g} & =  \sum\limits_{v=0}^{n-1} f(1^{\otimes v}\otimes g
\otimes 1^{\otimes n-v-1}) \nonumber\\ & -  (-1)^{ij+kl}
\sum\limits_{v=0}^{m-1} g(1^{\otimes v}\otimes f \otimes 1^{\otimes
m-v-1}) \quad\quad \in C^{n+m-1,i+j}_{k+l}(A,A) \nonumber
\end{align}

This is not the actual Lie bracket but the first step in our
construction. For degree and sign reasons we have to introduce
a shift map.

Let $S(A)$ be the bigraded module with $S(A)_{u}^{v} = A_{u}^{v+1}$,
and so the suspension map $S:A\to S(A)$ given by the identity map in each
bidgeree has internal bidegree $(0,-1)$. Given $f \in C^{n,i}_k(A,A)$,
then
\[
\sigma(f)=(-1)^{n+i+k-1} S\circ f \circ (S^{-1})^{\otimes n} \in
C^{n,i+n-1}_k(S(A),S(A)).
\]
Conversely, for $F \in C^{m,j}_l(S(A),S(A))$, we define
\[
\sigma^{-1}(F)=(-1)^{j+l+\binom{m}{2}}S^{-1}\circ F\circ S^{\otimes
m} \in C^{m,j+1-m}_l(A,A),
\]
so $\sigma^{-1}(\sigma(f))=f$. Particularly, for $m_{ij} \in
C^{j,2-(i+j)}_i(A,A)$, we have $\sigma(m_{ij})\in C^{j,1-i}_i(S(A),S(A))$.
Note that the notation $\sigma(f)$ does not mean applying a
shift functor to $f$.

We now define
\begin{align}
[f,g] & :=  \sigma^{-1}\bc{\sigma(f),\sigma(g)} \nonumber\\ & =
\sum\limits_{v=0}^{n-1} (-1)^{(n-1)(m-1)+v(m-1)+j(n-1)} f(1^{\otimes
v}\otimes g \otimes 1^{\otimes n-v-1}) \nonumber \\ & -
(-1)^{\left<f,g\right>} \sum\limits_{v=0}^{m-1}
(-1)^{(m-1)(n-1)+v(n-1)+i(m-1)} g(1^{\otimes v}\otimes f \otimes
1^{\otimes m-v-1}) \nonumber \\ & \qquad \in C^{n+m-1,i+j}_{k+l}(A,A)
\nonumber
\end{align}
for $f \in C^{n,i}_k(A,A)$ and $g \in C^{m,j}_l(A,A)$. Here,
$\left<f,g\right> := (n+i-1)(m+j-1)+kl$.
(See the Appendix for this computation.) It is easy to see that in the
case of bigraded modules concentrated in horizontal degree 0 this
specializes to the Lie algebra structure given in Section
\ref{liealgebrastructure}.
\medskip

As earlier, we use formal infinite sums of morphisms. These are now
bigraded and any such sum is actually finite in any given bidegree.

\begin{rrr}
It is also possible to work with a different definition of the shift
$\sigma$ on morphisms. Instead of our convention
\[
\sigma(f)=(-1)^{k+n+i-1} S \circ f\circ (S^{-1})^{\otimes n},
\]
it is also possible to work with
\[
\overline{\sigma}(f)=(-1)^{k+n+i-1} S \circ f\circ (S^{\otimes
n})^{-1}
\]
as in \cite[3.6]{Kel01} which differs from the above $\sigma$ by the
sign $(-1)^{\binom{n}{2}}$. Working with $\overline{\sigma}$ would
recover Keller's sign convention in the definition of
$A_\infty$-algebras and their morphisms, whereas our choice of
$\sigma$ recovers the signs of Lef\`evre-Hasegawa and Sagave.
\end{rrr}

\medskip
It is convenient to describe the above bracket in terms of a
composition product as in~\cite{Ger63}.

\begin{ddd} For $f\in C^{n,i}_k(A,A)$ and $g \in C^{m,j}_l(A,A)$ we
define the composition product $\circ$ by
\begin{align}
f \circ g  & =
\sum\limits_{v=0}^{n-1}\sigma^{-1}\Big((\sigma(f)(1^{\otimes
v}\otimes \sigma(g)\otimes 1^{\otimes n-v-1})\Big) \nonumber \\ & =
\sum\limits_{v=0}^{n-1} (-1)^{(m-1)(n-1)+v(m-1)+j(n-1)} f(1^{\otimes
v}\otimes g \otimes 1^{\otimes n-v-1}) \in C^{n+m-1,i+j}_{k+l}(A,A).
\nonumber
\end{align}
\end{ddd}

Hence, we have that
\[
[f,g]= f \circ g - (-1)^{\left<f,g\right>} g \circ f.
\]
\medskip

We will show that with this bracket $C^{*,*}_*(A,A)$ can be regarded as a bigraded Lie
algebra in the sense of the following definition.

\begin{ddd}
A bigraded $k$-module $X=\oplus X_i^j$ is a bigraded Lie algebra if there is a bracket operation
$[-,-]:X\otimes X\to X$ satisfying
\begin{itemize}
\item $[g,f]=-(-1)^{ab+kl}[f,g]$,
\item $(-1)^{ac+km}[[f,g],h]+(-1)^{ab+kl}[[g,h],f]+(-1)^{bc+lm}[[h,f],g]=0$,
\end{itemize}
for $f\in X^a_k$, $g\in X^b_l$, $h\in X^c_m$.
\end{ddd}

\pagebreak

\begin{ppp}\label{trigradedlie} The above bracket gives $C^{*,*}_*(A,A)$ the structure of
a bigraded Lie algebra for the bigrading where $f\in C^{n,i}_k$ is given bidegree
$(k, n+i-1)$; i.e.~for all $f,g,h\in C^{*,*}_*(A,A)$,
\begin{itemize}
\item $[g,f]=-(-1)^{\left<f,g\right>}[f,g]$,
\item $(-1)^{\left<f,h\right>} [[f,g],h] + (-1)^{\left<g,f\right>}
[[g,h],f] + (-1)^{\left<h,g\right>}[[h,f],g] = 0$.
\end{itemize}
\end{ppp}

\begin{proof}
The first point is immediate. For the graded Jacobi identity
 we will show that the composition product
$\circ$ makes $C^{*,*}_*(A,A)$ a bigraded pre-Lie ring in the sense that for $f \in
C^{n,i}_k(A,A)$, $g \in C^{m,j}_l(A,A)$ and $h \in C^{u,v}_w(A,A)$, we have
\begin{equation}\label{prelie}
(h \circ f) \circ g - (-1)^{\left<f,g\right>} (h \circ g) \circ f=
h\circ(f\circ g) - (-1)^{\left<f,g\right>} h \circ (g \circ f).
\end{equation}
We can then apply a direct computation analogous to the proof of
Theorem 1 of \cite{Ger63} which proves the claim. (For this, we note
that $\left<f\circ g,h\right> = \left<f,h\right>+\left<g,h\right>$.)

To prove the equation (\ref{prelie}), we note that
\[
f \circ g = \sigma^{-1}(\sigma(f)\odot\sigma(g))
\]
with
\[
F \odot G := \sum\limits_{r=1}^{n-1} F(1^{\otimes r} \otimes G
\otimes 1^{\otimes n-r-1}).
\]
This is going to simplify the signs in (\ref{prelie}) considerably
since this equation is equivalent to
\begin{equation}\label{shiftprelie}
(H \odot F) \odot G - (-1)^{\left<f,g\right>} (H\odot G)\odot F = H
\odot (F \odot G) - (-1)^{\left<f,g\right>} H \odot (G \odot F)
\end{equation}
for $F = \sigma(f)$, $G=\sigma(g)$ and $H=\sigma(h)$. We have
\begin{align}
(H \odot F) \odot G
    & =  \big( \sum\limits_{r=0}^{u-1} H(1^{\otimes
            r} \otimes F \otimes 1^{\otimes u-r-1}) \big) \odot G \nonumber \\
    & = (-1)^{\left<f,g\right>} \sum\limits_{r=0}^{u-1}
            \sum\limits_{a+b=r-1} H(1^{\otimes a} \otimes G \otimes 1^{\otimes
            b} \otimes F \otimes 1^{\otimes u-r-1}) \nonumber\\
    & + \sum\limits_{r=0}^{u-1} \sum\limits_{s=0}^{n-1} H(1^{\otimes r}
            \otimes
            F(1^{\otimes s} \otimes G \otimes 1^{\otimes n-s-1}) \otimes 1^{\otimes u-r-1}) \nonumber \\
    & +  \sum\limits_{r=0}^{u-1} \sum\limits_{a+b=u-r-2} H(1^{\otimes r}
            \otimes F \otimes 1^{\otimes a} \otimes G \otimes 1^{\otimes b}).
\nonumber
\end{align}
Note that the sign $(-1)^{\left<f,g\right>}$ in the first summand
arises from the Koszul sign rule for interchanging $F$ and $G$. Using
this, we can read off the equation (\ref{shiftprelie}), from which
(\ref{prelie}) follows.
\end{proof}

Now we would like to describe derived $A_\infty$-structures in terms
of this Lie algebra structure, but first we have to introduce
another operation which alters signs.

\begin{ddd}
For $f \in C^{n,i}_k(A,A)$ define $f^\#=(-1)^k f \in C^{n,i}_k(A,A)$.
\end{ddd}

This operation satisfies
\begin{itemize}
\item $(f^\#)^\#=f$,
\item $(f\circ g)^\# = f^\# \circ g^\#$,
\item $[f,g]^\#=[f^\#,g^\#]$.
\end{itemize}

\begin{ppp}\label{dainfzero}
Let $A$ be a bigraded $k$-module without 2-torsion with given map $\eta:k\to A$.
Let $m=\sum\limits_{\substack{i \ge 0, \, j \ge 1}} m_{ij}$ with
$m_{ij} \in C^{j,2-(i+j)}_i(A,A)$ satisfying the unit conditions of
Definition~\ref{derivedainfstructure}. 

\pagebreak

Then the following are equivalent
\begin{itemize}
\item $m$ is a derived $A_\infty$-structure on A,
\item $m\circ m^\# = 0$,
\item $[m,m^\#]=0$.
\end{itemize}
\end{ppp}

\begin{proof}
The equivalence of the first two points follows immediately from the
definitions. For the equivalence of the last two points let us
consider the part $[m,m^\#]_u$ of $[m,m^\#]$ that lies in horizontal
degree $u$. We have
\begin{align}
[m,m^\#]_u & =  \sum\limits_{\substack{u=i+p}} \big( m_{ij} \circ
m_{pq}^\# - (-1)^{ip+(i-1)(p-1)} m_{ij}^\# \circ m_{pq} \big)
\nonumber \\
& =  \sum\limits_{\substack{u=i+p}} \big((-1)^p m_{ij} \circ m_{pq}
- (-1)^{u+1+i} m_{ij} \circ m_{pq} \big). \nonumber
\end{align}
We are going to distinguish between the cases $u$ even and $u$ odd.
For even $u=i+p$, the sum splits into the cases where either both
$i$ and $p$ are even or both $i$ and $p$ are odd. In either case, we
can read off that
\[
[m,m^\#]_u=2(m \circ m^\#)_u.
\]
The case of $u$ odd follows similarly.
\end{proof}

\section{Hochschild cohomology and uniqueness of derived $A_\infty$-algebras}
\label{sec:deruniq}
\subsection{Hochschild cohomology of $dA_\infty$-algebras}

We would like to define a notion of Hochschild cohomology for
$dA_\infty$-algebras that extends the classical, non-derived case.
However, this is not as straightforward as before. In the classical
case of an $A_\infty$-algebra $A$ with $A_\infty$-structure $m$, we
could define a differential on $C^{*,*}(A,A)$ via $D=[m,-]$. This
satisfies $D \circ D = [m,[m,-]]=0$ since $[m,m]=0$. But in the
derived case the signs are slightly more complicated which means we
can only guarantee $[m,m^\#]=0$. We can still define Hochschild
cohomology for a certain class of $dA_\infty$-algebra which includes
the cases we are interested in.

\medskip
\begin{ddd}
Let $m=\sum\limits_{i \ge 0, j\ge 1} m_{ij}$ be a
$dA_\infty$-structure. Then we denote the horizontal even degree
part by $m_{even}$ and the horizontal odd degree part by $m_{odd}$,
i.e.,
\[
m_{even}=\sum\limits_{i \,\mbox{\scriptsize{even}}} m_{ij}
\quad\quad\mbox{and} \quad\quad m_{odd}=\sum\limits_{i
\,\mbox{\scriptsize{odd}}} m_{ij}.
\]
\end{ddd}

\begin{rrr}
Since $m$ is a $dA_\infty$-structure, by Lemma \ref{dainfzero} we
have $(m_{even}+m_{odd})\circ(m_{even}-m_{odd})=0,$ which splits as
\[
m_{even}\circ m_{even} = m_{odd} \circ m_{odd}
\quad\quad\mbox{and}\quad\quad m_{even}\circ m_{odd}=m_{odd}\circ
m_{even}.
\]
\end{rrr}

\begin{ddd}
We call a derived $A_\infty$-structure $m$ \emph{orthogonal} if
$$m_{even}\circ m_{even}=0 \quad\quad\mbox{or, equivalently,}\quad\quad m_{odd}\circ
m_{odd}=0.$$
\end{ddd}

\begin{ex}
Bidgas are orthogonal since they have $m_{odd}=m_{11}$ and $m_{11}\circ m_{11}=0$.
\end{ex}

\begin{lem}\label{derivedzerobracket}
Let $A$ be a bigraded $k$-module without 2-torsion and let $m=\sum\limits_{i,j}
m_{ij}$ be an orthogonal derived $A_\infty$-structure on $A$. Define
\[
D: C^{*,*}_*(A,A) \longrightarrow C^{*,*}_*(A,A)
\]
via
\[
D(f)=[m_{even},f^\#] + [m_{odd},f] = (-1)^k [m_{even},f] +
[m_{odd},f] \quad\quad\mbox{for} \quad f\in C^{n,i}_k(A,A).
\]
Then $D$ satisfies $D\circ D = 0$. Also, $D$ raises the total degree
by 1, so $D$ is a differential on $C^{*,*}_*(A,A)$.
\end{lem}

\begin{proof}
The map $D$ raises degree by 1 since $m$ has total degree 2. Let us
look at $D(D(f))$. Assume that $f$ has horizontal internal degree
$k$. Then for even $p$ the horizontal degree of $[m_{pq},f]$ has the
same parity as $k$ whereas for odd $p$ the horizontal degree of
$[m_{pq},f]$ has the parity of $k+1$. This means that
\[
[m_{even},f]^\#=(-1)^k[m_{even},f] \quad\quad\mbox{and}\quad\quad
[m_{odd},f]^\#=(-1)^{k+1}[m_{odd},f].
\]
Thus, we obtain
\[
D((-1)^k[m_{even},f])=(-1)^k\big((-1)^k[m_{even},[m_{even},f]]+[m_{odd},[m_{even},f]]
\big)
\]
and
\[
D([m_{odd},f])=(-1)^{k+1}[m_{even},[m_{odd},f]]+[m_{odd},[m_{odd},f]]
\]
which together give us
\[
D(D(f))= [m_{even},[m_{even},f]]+(-1)^k[m_{odd},[m_{even},f]] +
(-1)^{k+1}[m_{even},[m_{odd},f]]+[m_{odd},[m_{odd},f]].
\]
Since $m$ is assumed to be orthogonal, we can directly compute that
\[
[m_{even},[m_{even},f]]=0=[m_{odd},[m_{odd},f]].
\]
From the graded Jacobi identity established in Proposition
\ref{trigradedlie} we conclude that
\[
[m_{odd},[m_{even},f]]=[m_{even},[m_{odd},f]].
\]
Putting this together, we can read off the desired equation $D\circ D=0$.
\end{proof}

\begin{ddd}
Let $A$ be an orthogonal $dA_\infty$-algebra with orthogonal
$dA_\infty$-structure $m$. Then the Hochschild cohomology of $A$ as
a $dA_\infty$-algebra is defined as
\[
\hh^{*}(A,A) := H^*\left(\prod\limits_{i,j} C^{i,*-i-j}_j(A,A),D\right).
\]
\end{ddd}

The grading in the above definition of Hochschild cohomology denotes
the total degree.

\begin{rrr}
If $A$ has $dA_\infty$-structure $m=m_{11}+m_{02}$ (i.e.~$A$ is a
bidga with trivial vertical differential), then this definition
specializes to Sagave's definition~\cite[Section 5]{Sag10} of
Hochschild cohomology of bidgas with trivial vertical differential.
\end{rrr}

In this very special case of a bidga with trivial vertical differential,
one grading is preserved by both $m_{11}$ and $m_{02}$ so that 
we have bigraded Hochschild cohomology groups:
\[
\hh^t(A,A) = \prod\limits_{s \geq 0} \hh^{s,t-s}(A,A),
\]
where $\hh^{s,r}(A,A)=H^s(\prod_n C^{n,r}_{*-n}(A,A),D)$.
We denote the Hochschild cohomology in this special case by
$\hh^{*,*}_{bidga}(A,A)$.

\subsection{Uniqueness of derived $dA_\infty$-algebras}

The overall goal of this section is to establish a uniqueness result
analogous to Kadeishvili's (Theorem \ref{Kadeishviliuniqueness}) for the possibility of
extending an existing $dA_\infty$-structure on a minimal model. A
minimal model of a differential graded algebra has an underlying
bidga with zero vertical differential. Let $\mu=m_{02}$ denote the
multiplication of this bidga and $\partial=m_{11}$ the horizontal
differential.

The first step is to look into how to perturb an existing
$dA_\infty$-structure by certain elements $b$ of total degree 1.
\begin{ddd}
Let $A$ be a bidga with multiplication
$m_{02}=\mu$, horizontal differential $m_{11}=\partial$ and vertical
differential $m_{01}=0$. Then
$$a = \sum\limits_{i\ge0, \,j\ge1} a_{ij}, \quad \quad a_{ij} \in C^{j,2-(i+j)}_i(A,A),
\quad i+j \ge 3,$$ is a {\it twisting cochain} if $\partial + \mu +
a$ is a $dA_\infty$-structure.
\end{ddd}

\begin{rrr}
Note that by Proposition~\ref{dainfzero} $a$ is a twisting cochain if and only if we have
\[
[\partial + \mu + a, \partial^\# + \mu^\# + a^\#]=0.
\]
Letting $D$ be the differential corresponding to the orthogonal
$dA_\infty$-structure $m=\partial +\mu$, this is equivalent to the
\emph{derived Maurer-Cartan formula}
\begin{equation}\label{derivedmaurercartan}
2D(a)=-[a,a^\#] + 4[\partial,a_{odd}],
\end{equation}
as can be verified quickly by splitting $a$ into even and odd
horizontal degree parts and
using that $[\partial+\mu, \partial+\mu]=0$. Hence, an element $a = \sum\limits_{i,j}
a_{ij}, \quad a_{ij} \in C^{j,2-(i+j)}_i(A,A), \quad i+j \ge 3$ is a
twisting cochain if and only if $a$ satisfies the above derived Maurer-Cartan
formula.
\end{rrr}

\begin{lem}\label{derivedperturb} Let $A$ be a bidga with multiplication
$m_{02}=\mu$, horizontal differential $m_{11}=\partial$ and vertical
differential $m_{01}=0$. Let
$$a = \sum_{i,j} a_{ij}, \quad \quad a_{ij} \in C^{j,2-(i+j)}_i(A,A),
\quad i+j \ge 3,$$ be a twisting cochain. Let either
\begin{description}
\item[(A)] $b \in C^{n-1,2-(n+k)}_k(A,A)$, for $k+n \ge 3$, with
$[\partial,b]=0$
\end{description}
or
\begin{description}
\item[(B)] $b \in C^{n,2-(n+k)}_{k-1}(A,A)$, for $k+n \ge 3$, with
$[\mu,b]=0$.
\end{description}

Then there is a twisting cochain $\overline{a}$ satisfying
\begin{itemize}
\item the $dA_\infty$-structures $\partial + \mu + a$ and $\overline{m} = \partial + \mu + \overline{a}$ are
$E_2$-equivalent,
\item $\overline{a}_{uv} = a_{uv}$ for $u<k$ or $v<n-1$ or $(u,v)=(k,n-1)$ in case {\bf (A)} and for $u<k-1$ or $v<n$
or $(u,v)=(k-1,n)$ in case {\bf (B)},
\item $\overline{a}_{kn} = a_{kn} - [\mu,b]$ in
case {\bf (A)},
\item $\overline{a}_{kn}= a_{kn} - [\partial,b] $ in case {\bf (B)}.
\end{itemize}
\end{lem}

\begin{proof}
This is a lengthy but direct computation using the definition of a 
morphism of $dA_\infty$-algebras. The twisting cochain $\overline{a}$ is going to be determined
by $\partial+\mu+a$ being $E_2$-equivalent to $\partial+\mu+\overline{a}$ via
the equivalence $id+b$.
We will only do case {\bf (A)} explicitly since the
other case can be read off the proof of this one.

Let $f:= id + b$. We consider what it means for there to be a
$dA_\infty$-structure $\overline{m}= \partial +\mu + \overline{a}$
on $A$ such that $f: (A,m) \longrightarrow (A,\overline{m})$ is a
morphism of $dA_\infty$-structures, i.e.~the
equation~(\ref{dmorphequation}) in Definition
\ref{derivedainfmorphism} is satisfied. Using $f_{01}= id$,
$f_{k,n-1}=b$ and $f_{ij}=0$ in all other degrees as well as $m =
\mu + a$ and $\overline{m} = \mu + \overline{a}$, we write down
(\ref{dmorphequation}). The left-hand side of (\ref{dmorphequation})
is only nonzero for $(i,j)=(0,1)$ and $(i,j)=(k,n-1)$. Thus, we
obtain
\[
(-1)^u m_{uv} + \sum\limits_{r=0}^{n-2}
(-1)^{r(v-n)+(n-r)+(u-k)(n-1)} b(1^{\otimes r} \otimes m_{u-k,
v+2-n} \otimes 1^{\otimes n-2-r}).
\]
The sum can only be nonzero if $u \ge k$ and $v \ge n-1$ and
$(u,v)\neq(k,n-1)$. In the special case $(u,v) = (k,n)$ we get
\[
(-1)^k a_{kn} + \sum\limits_{r=0}^{n-2} (-1)^{n-r} b(1^{\otimes r}
\otimes \mu \otimes 1^{\otimes n-2-r}).
\]
For $(u,v)=(k+1,n-1)$, the result is
\[
(-1)^{k+1} a_{k+1,n-1} - \sum\limits_{r=0}^{n-2}  b(1^{\otimes r}
\otimes \partial \otimes 1^{\otimes n-2-r}).
\]

On the right-hand side of (\ref{dmorphequation}) we have
\begin{equation}\label{rhs}
(-1)^{u} \overline{m}_{uv} + \sum\limits_{\substack{ u= i + p_1 +\cdots + p_j \\
 v= q_1+\cdots + q_j}} (-1)^{\epsilon} \overline{m}_{ij}
(f_{p_1 q_1} \otimes \cdots \otimes f_{p_j q_j})
\end{equation}
where at least one of the $f_{p_r q_r}$ in the sum has to be
$f_{k,n-1}=b$ and $\epsilon$ is as in Definition
\ref{derivedainfmorphism}. The following four special cases are to
be considered. First, we note that, since we have $\overline{m}_{01}=0$, the sum is zero for
$(u,v)=(k,n-1)$. For $(u,v) = (k,n)$, we obtain
\[
(-1)^{k} \overline{a}_{kn} + (-1)^{n} \mu(1 \otimes b) + \mu(b
\otimes 1),
\]
for $(u,v)=(k+1,n-1)$ we have
\[
(-1)^{k+1}\overline{a}_{k+1,n-1} + (-1)^{k+1}\partial(b)
\]
and for $(u,v) = (2k,2n-2)$ the result is
\[
 \overline{a}_{2k, 2n-2} + (-1)^{nk} \mu(b\otimes b) + \sum\limits_{r=0}^{n-1}
(-1)^{\epsilon} \overline{a}_{k,n} (1^{r} \otimes b \otimes
1^{n-1-r}).
\]
In all other cases each summand appearing in the sum in (\ref{rhs}) has $i+j \ge 3$.
Further, the sum in (\ref{rhs}) can only be nonzero for $u \ge i+k$
and $v \ge (n-1) + (j-1)$.

Now recall that
\[
[\partial,b]=\partial(b)-(-1)^k \sum\limits_{r=0}^{n-2} b(1^{\otimes
r} \otimes
\partial \otimes 1^{\otimes n-2-r})
\]
and
\[
[\mu,b]= (-1)^{n+k}\Big(\mu(1\otimes b) + (-1)^{n}\mu(b\otimes 1) +
 \sum\limits_{r=0}^{n-2} (-1)^{r+1} b(1^{\otimes r} \otimes\mu\otimes
1^{\otimes n-2-r}) \Big).
\]
Further, note that we have assumed that $[\partial,b]=0$.

Putting all this together, we can read off that for $(u,v)$ with
either $u<k$ or $v<n-1$ and for $(u,v)=(k,n-1)$, we have
\[
\overline{a}_{uv} = a_{uv}.
\]
For $(u,v)$=$(k,n)$, we get
\begin{align}
\overline{a}_{kn} & =  a_{kn} - (-1)^{k}\Big( \mu(b\otimes 1)
+(-1)^n \mu(1\otimes b) + (-1)^{n-1} \sum\limits_{r=0}^{n-2}
(-1)^{r}
b(1^{\otimes r} \otimes\mu\otimes 1^{\otimes n-2-r}) \Big)\nonumber\\
&= a_{kn} - [\mu,b]; \nonumber
\end{align}
for $(u,v)=(k+1,n-1)$ we have
\begin{align}
\overline{a}_{k+1,n-1} & = a_{k+1,n-1} + (-1)^k
\sum\limits_{r=0}^{n-2} b(1^{\otimes r} \otimes
\partial \otimes 1^{\otimes n-2-r}) - \partial(b) \nonumber\\ & = a_{k+1,n-1}
-
[\partial,b] = a_{k+1,n-1}; \nonumber
\end{align}
for $(u,v)=(2k,2n-2)$ we have
\begin{align}
\overline{a}_{2k,2n-2} & = a_{2k,2n-2} +
\sum\limits_{r=0}^{n-2} (-1)^{m+n-r+k(n-1)}
b(1^{\otimes r} \otimes
a_kn \otimes 1^{\otimes n-2-r}) -
(-1)^{nk}\mu(b\otimes b)\nonumber \\
&\qquad+\sum_{r=0}^{n-1} (-1)^{\epsilon}\overline{a}_{kn}(1^{\otimes r}\otimes b\otimes 1^{\otimes n-1-r}). \nonumber
\end{align}
Finally for $(u,v)\neq (k,n), (k+1, n-1)$ or $(2k,2n-2)$ with $u \ge k$ and $v \ge n-1$, we have
\begin{align}
    \overline{a}_{uv} & = a_{uv} + (-1)^{u} \sum\limits_{r=0}^{n-2}
                            (-1)^{r(v-n)+(n-r)+(u-k)(n-1)} b(1^{\otimes r} \otimes m_{u-k,
                            v+2-n} \otimes 1^{\otimes n-2-r}) \nonumber\\
                        & - (-1)^{u}\sum\limits_{\substack{ u= i + p_1 +\cdots + p_j \\ v= q_1+\cdots+ q_j \\
                            \mbox{\scriptsize{at least one $q_j\neq 1$}}}} (-1)^{\epsilon}
                            \overline{a}_{ij} (f_{p_1 q_1} \otimes \cdots \otimes f_{p_j q_j}).
                            \nonumber
\end{align}
Note that the second sum in the last equation can only be nonzero if
$i+j \ge 3$, $u\ge k+i$ and $v \ge (n-1)+(j-1)$. Also, for fixed
$(u,v)$, the right-hand side of the last equation only uses
$\overline{a}_{pq}$ with $p<u$ and $q<v$. The same thing happens in the case
$(u,v)=(2k, 2n-2)$.
This proves that the
$\overline{a}$ in the statement of our lemma can be constructed
inductively.

One can then check degreewise that $\overline{m}=
\partial + \mu + \overline{a}$ defines a $dA_\infty$-structure by showing
that $[\overline{m},\overline{m}^\#]=0$. The morphism $f$ is an
$E_2$-equivalence since $f_{01}=id$.
\end{proof}

\begin{rrr}
Note that in the situation of the above lemma, in both cases we have in particular that $\overline{a}_{uv}=a_{uv}$
whenever $u+v<k+n$.
\end{rrr}

We can now formulate a derived version of Kadeishvili's uniqueness
theorem.

\begin{theorem}
\label{derivedHHcriterion} Let $A$ be a bidga with multiplication
$m_{02}=\mu$, horizontal differential $m_{11}=\partial$ and vertical
differential $m_{01}=0$. If $\hh^{r,2-r}_{bidga}(A,A)=0$ for
$r \ge 3$, then every $dA_\infty$-structure on $A$ with
$m_{02}=\mu$, $m_{11}=\partial$ and $m_{01}=0$ is $E_2$-equivalent
to the trivial one, i.e.~the one with $m_{02}=\mu$,
$m_{11}=\partial$ and $m_{ij}=0$ for $(i,j) \neq (0,2)$ or $(1,1)$.
\end{theorem}

\begin{proof}
Let $m = \partial + \mu + a$ be a $dA_\infty$-structure on $A$ with
\[
a=\sum\limits_{k+n \ge 3} a_{kn}, \quad\quad a_{kn} \in
C^{n,2-(k+n)}_k(A,A).
\]
We want to show that $m$ is $E_2$-equivalent to the
$dA_\infty$-structure $\partial + \mu$.

\medskip
We now fix $t \ge 3$ and show that $m$ is equivalent to a
$dA_\infty$-structure with $a_{kn}=0$ for $k+n=t$. We show this by
induction on $k$. Assuming that $a_{ij}=0$ for $i+j=t$ and $i<k$, we
will show that $m$ is equivalent to a $dA_\infty$-structure with
$\overline{m}=\partial + \mu + \overline{a}$ with
$\overline{a}_{kn}=0$ and $\overline{a}_{ij}=a_{ij}=0$ for $i+j=t$,
$i<k$ and $i+j <t$.

\medskip
Because $m$ is a $dA_\infty$-structure, by Lemma
\ref{derivedzerobracket} we have $[\partial + \mu + a,\partial^\# +
\mu^\# + a^\#]$=0. Since $A$ is also a bidga, we have $[\partial +
\mu,
\partial^\# + \mu^\#]=0$. Hence, $a$ is a twisting cochain satisfying the Maurer-Cartan formula
\[
2D(a)=-[a,a^\#]+4[\partial,a_{odd}]
\]
as explained in (\ref{derivedmaurercartan}). Further, we have
\[
D(-) = [\mu,(-)^\#] + [\partial,-]
\]
with
\[
[\mu,(-)^\#]: C^{*,*}_*(A,A) \longrightarrow C^{*+1,*}_*(A,A)
\quad\quad\mbox{and} \quad\quad  [\partial,-]: C^{*,*}_*(A,A)
\longrightarrow C^{*,*}_{*+1}(A,A),
\]
so $[\mu, a_{kn}^\#]$ lives in the tridegree $(n+1,k,2-(k+n))$-part
of $D(a)$ and $[\partial,a_{kn}]$ lives in tridegree
$(n,k+1,2-(k+n))$. However, on the other side of (\ref{derivedmaurercartan}) the tridegree
$(n+1,k,2-(k+n))$-part as well as the $(n,k+1,2-(k+n))$-part of
$[a,a^\#]$ is zero since $[a,a^\#]$ can only be nonzero in degrees
$(u,v,w)$ with $u+v \ge 5$ whereas $n+1+k = 4$. Here we are adopting the
convention for tridegrees that an element in $C^{n,i}_k(A,A)$ has
tridegree $(n,k,i)$.

Thus according to
the Maurer-Cartan formula, $D(a_{kn})$ lives in
$2[\partial,a_{odd}]$. This information splits into the equations
\begin{equation}
[\mu,a_{kn}^\#]= \epsilon_1 2 [\partial,a_{k-1,n+1}], \quad\quad
\epsilon_1\in\{0,1\} \quad\quad\mbox{and}\quad\quad
[\partial,a_{kn}]= \epsilon_2 2 [\partial,a_{kn}], \quad\quad
\epsilon_2\in\{0,1\} \nonumber
\end{equation}
where $\epsilon_2=0$ for $k$ even by definition (since the right
hand side is supposed to be a summand of $2 [\partial,a_{odd}]$).
Thus, we can also conclude that $[\mu,a_{kn}^\#]=0$ since our
induction assumption gives $a_{k-1,n+1}=0$.

\medskip
For $k$ odd, we are left with $[\partial,a_{kn}]= \epsilon_2 2
[\partial,a_{kn}],\, \epsilon_2\in\{0,1\}$, from which we can
immediately read off that $[\partial,a_{kn}]=0$.

\medskip
Hence, in any case $D(a_{kn})=0$ and $a_{kn}$ is a cocycle in
$C^{n,2-(n+k)}_k(A,A)$, so $$[a_{kn}]\in \hh^{k+n,2-k-n}_{bidga}(A,A).$$ 
However, $\hh^{k+n,2-k-n}_{bidga}(A,A)$ is zero by
assumption, so there must be a $b$ in total degree 1 with $D(b) =
a_{kn}$.

\medskip So, analogously to the proof of Theorem \ref{th:HHcriterion}, there is
a $b_1 \in C^{n-1,2-(k+n)}_{k}(A,A)$ with $[\partial,b_1]=0$ and
$[\mu,b_1]=a_{kn}$ and $b_2 \in C^{n,2-(k+n)}_{k-1}(A,A)$ with
$[\mu,b_2]=0$ and $[\partial,b_2]=a_{kn}$.

\medskip Applying Lemma \ref{derivedperturb} to $b_1$, there is a
$dA_\infty$-structure $\overline{m}= \partial + \mu +
\overline{a}_{ij}$
 with $\overline{a}_{ij} \in C^{j,2-(i+j)}_i(A,A)$, $i+j
\ge 3$ such that $\overline{m}$ is $E_2$-equivalent to $m$,
$\overline{a}_{kn} = a_{kn}-[\mu, b_1]=0$ and $\overline{a}_{ij}=a_{ij}$ for $i+j <t$
and $i+j=t$, $i<k$, which proves our claim.
\end{proof}

\begin{ex}
In \cite[Proposition 4.2]{DugShi09}, Dugger and Shipley consider the
dga
\begin{eqnarray}
A=\mathbb{Z}\left<e,x,y\right>\big/(e^2=0, ex+xe=x^2, xy=yx=1),
\nonumber \\ \partial(e)=p, \partial(x)=0, \partial(y)=0, \quad\quad
|e|=|x|=1, |y|=-1. \nonumber
\end{eqnarray}
This is a dga over $\mathbb{Z}$ which has homology $H_n(A)=\mathbb{Z}/p$ in every
degree $n$. (Note that Dugger and Shipley use homological grading.)
They then prove in Theorem 4.5 that $A$ is not formal.

In \cite{Sag09} Sagave gives a projective $E_1$-resolution $B$ of
$A$. He then  constructs the first degrees of a minimal model
structure on the induced termwise projective resolution
$P=H^*_{ver}(B)$ and shows that this gives a nontrivial class in
$\hh^{3,-1}_{bidga}(P,P)$.
\end{ex}

Theorem~\ref{derivedHHcriterion} will be used in the next section to give a sufficient
criterion for the existence of a unique dga realising a fixed
homology algebra over a ground ring rather than a ground field. To
prove this derived analogue of Corollary~\ref{uniquedga}, we first
have to investigate the behaviour of Hochschild cohomology of
degreewise projective resolutions under $E_2$-equivalence.

\section{Invariance under $E_2$-equivalence and intrinsic formality}
\label{sec:e2inv}

In order to establish our uniqueness criterion we need an invariance result
for Hochschild cohomology under $E_2$-equivalence. To prove this
we will need to define Hochshild cohomology with coefficients.
We will carry this out here only for the special case we need.
In future work we hope to study the general case, but this
would take us too far afield here.

Thus we will concentrate on the case of relevance to us, namely
bidgas with $m_{01}=0$. Invariance under
$E_2$-equivalence in this situation is also discussed in~\cite[Section 5]{Sag10}. We begin
by spelling out concretely what a bidga with $m_{01}=0$ is.
\smallskip

A bidga with $m_{01}=0$ is a bigraded module $A_i^j$ equipped with maps
$m_{11}:A_i^j\to A_{i-1}^j$ and $m_{02}:(A\otimes A)_i^j\to A_i^j$ with
relations which specify that $m_{02}$ is associative, $m_{11}$ is a
differential and $m_{11}$ is a derivation with respect to $m_{02}$.
These relations come from the cases $(u=0, v=3)$, $(u=2, v=1)$ and $(u=1, v=2)$
respectively of the defining relations; all other relations are trivial.
Notice that this is just a dga with an extra grading.

It is straightforward to see what a module over such a thing should be; it is just
a dg module with an extra grading.

\begin{ddd}
Let $A$ be a bidga with $m_{01}=0$. A \emph{left $A$-module} $M$ is a bigraded module $\{M_i^j\}$
over the ground ring equipped with a horizontal differential $\overline{m}_{11}:M_i^j\to
M_{i-1}^j$ and an associative action $\overline{m}_{02}^l:(A\otimes M)_i^j\to M_i^j$
such that the diagram
    $$
    \xymatrix{
    A\otimes M \ar[r]^{\overline{m}_{02}^l}
    \ar[d]_{m_{11}\otimes 1+1\otimes \overline{m}_{11}} & M \ar[d]^{\overline{m}_{11}}\\
    A\otimes M \ar[r]_{\overline{m}_{02}^l} & M
    }
    $$
commutes.
\end{ddd}

A right $A$-module is defined in the obvious way, with a right action
map $\overline{m}_{02}^r:M\otimes A\to M$. And an $A$-bimodule is simultaneously
a left and right $A$-module with the obvious compatibility condition on the
left and right actions.

Notice that a morphism of bidgas $A\to A'$ between bidgas with $m_{01}=0$
makes $A'$ into an $A$-bimodule.

\smallskip

Let us also spell out what an $E_2$-equivalence $f:A\to A'$ between bidgas with $m_{01}=0$
is. This is just a morphism $f:A\to A'$ inducing an isomorphism on horizontal homology.
(So we can think of such an $f$ as a quasi-isomorphism if we think of $A$ and $A'$
as complexes with respect to horizontal differentials.)

\medskip

Now let $A$ be a bidga with $m_{01}=0$ and let $M$ be an $A$-bimodule.
Let
    $$
    C_k^{n,i}(A,M)=\prod_{u,v} \Hom\left((A^{\otimes n})_u^v,M_{u-k}^{v+i}\right)
    $$
and for $f\in C_k^{n,i}(A,M)$ define
    \begin{align*}
    Df=&(-1)^{k+n+i-1} \overline{m}_{02}^r(f\otimes 1) +(-1)^{k+i} \overline{m}_{02}^l(1\otimes f)\\
    &\quad +(-1)^{k+n+i} f\circ m_{02} +\overline{m}_{11}\circ f+(-1)^{k+1}f\circ m_{11}.
    \end{align*}

Then $D$ is a differential, allowing us to make the following definition.

\begin{ddd}
For $A$ a bidga with $m_{01}=0$ and $M$ an $A$-bimodule the \emph{Hochschild
cohomology of $A$ with coefficients in $M$} is defined by
    $$
    \hh^{s,r}_{bidga}(A,M)=H^s\left(\prod_n C^{n,r}_{*-n}(A,M),D\right).
    $$
\end{ddd}

This is a covariant functor of $M$ and a contravariant functor of $A$.
In the case where $M=A$, regarded as a bimodule over itself, this agrees
with the earlier definition of $\hh^{*,*}_{bidga}(A,A)$. Indeed the formula above for
the differential $D$ just becomes $Df=(-1)^k[m_{02},f]+[m_{11},f]$.

\begin{ppp}
\label{E2inv} Let $(A,m)$ and $(A',\overline{m})$ be bidgas with
$m_{01}=\overline{m}_{01}=0$ and which are degreewise projective
over $k$. Let $f:A\to A'$ be an $E_2$-equivalence. Then
$f$ induces an isomorphism of Hochschild cohomology groups
    $$\hh^{*,*}_{bidga}(A,A)\cong \hh^{*,*}_{bidga}(A',A').$$
\end{ppp}

\begin{proof}
For each $i$,
we can interpret the Hochschild cohomology $HH^{*,i}_{bidga}(A,M)$  as
the cohomology of a (right half-plane) bicomplex.  
This works very similarly to the case
of Hochschild cohomology of a dga discussed earlier. One differential, say $D_1$,
is given by the $m_{11}$ part of the formula for $D$ and the other, say $D_2$,
by the $m_{02}$ part.

Now consider $A$ and $M$ as complexes with respect to their differentials
$m_{11}$ and $\overline{m}_{11}$ (with an extra grading). The differential $D_1$ 
on $\Hom_*^i(A^{\otimes p}, M)$
is the induced
differential via the tensor product and $\Hom$ functors of complexes.
For bounded below and degreewise
projective complexes the ordinary $\Hom$ and tensor product functors agree with the
derived versions and are therefore quasi-isomorphism invariant.

Thus the morphism $f:A\to A'$ induces column-wise quasi-isomorphisms
of bicomplexes $C^{*,i}_*(A,A)\to C^{*,i}_*(A,A')$ and
$C^{*,i}_*(A',A')\to C^{*,i}_*(A,A')$. 
It follows that the induced maps of total complexes are quasi-isomorphisms
and $\hh^{*,*}_{bidga}(A,A)\cong \hh^{*,*}_{bidga}(A',A')$.

\end{proof}

Now we are in a position to give our criterion for intrinsic formality.

\begin{theorem}
\label{th:derivedformality} Let $A$ be a dga and $E$ its minimal
model with $dA_\infty$-structure $m$. By $\tilde{E}$, we denote the
underlying bidga of $E$, i.e. $\tilde{E}=E$ as $k$-modules together
with $dA_\infty$-structure $\tilde{m}=m_{11}+m_{02}$. If
\[
\hh^{m,2-m}_{bidga}(\tilde{E},\tilde{E})=0
\quad\quad\mbox{for}\quad\quad m \ge 3,
\]
then $A$ is intrinsically formal.
\end{theorem}

\begin{proof}
Applying Theorem~\ref{derivedHHcriterion} to $\tilde{E}$, we obtain
that every $dA_\infty$ structure on $\tilde{E}$ is $E_2$-equivalent
to the trivial one. By definition of minimal model, $A$ is
$E_2$-equivalent to $E$. Thus $A$ is $E_2$-equivalent to
$(\tilde{E}, \text{triv})$. Again by definition of minimal model,
$(\tilde{E}, \text{triv})$ is $E_2$-equivalent to $(H^*(A),
\text{triv})$. Thus we have an $E_2$-equivalence between $A$ and
$(H^*(A), \text{triv})$ and since these are both dgas an
$E_2$-equivalence is a quasi-isomorphism. So $A$ is formal.

Now let $A'$ be a dga with $H^*(A)\cong H^*(A')$ as associative algebras,
let $E'$ be a minimal model of $A'$ and let $\tilde{E'}$ be its underlying
bidga. We have $E_2$-equivalences
    $$
    \tilde{E'}\simeq (H^*(A'),\text{triv})\simeq (H^*(A), \text{triv})\simeq \tilde{E}.
    $$
Thus $\tilde{E'}$ and $\tilde{E}$ are $E_2$-equivalent bidgas. By
definition of minimal model they are degreewise projective and have
$m_{01}=0$. Applying Proposition~\ref{E2inv} gives
$\hh^{m,2-m}_{bidga}(\tilde{E'}, \tilde{E'})\cong
\hh^{m,2-m}_{bidga}(\tilde{E}, \tilde{E})$. So the Hochschild
cohomology of $\tilde{E'}$ is zero in the relevant range and the
argument of the preceeding paragraph shows that $A'$ is also formal.

Since $A$ and $A'$ are both formal, the hypothesis $H^*(A)\cong H^*(A')$ means
they are quasi-isomorphic.
\end{proof}

\section{Uniqueness of classical $A_\infty$-structures}
\label{sec:unique}

In this section $k$ is still a commutative ground ring without 2-torsion unless stated otherwise.
We use Hochschild cohomology of differential graded algebras to give a uniqueness 
criterion for extending the differential and multiplication of a fixed dga to an $A_\infty$-structure. 
In the case of a trivial differential this recovers Kadeishvili's classical 
Theorem~\ref{Kadeishviliuniqueness}. We then apply this to an example in homotopy theory. 

\bigskip
Fix a differential graded algebra $A$ with differential
$m_1=\partial$ and multiplication $m_2=\mu$. We would like to
consider the set of all $A_\infty$-structures on $A$ (up to
quasi-isomorphism) that extend the differential graded algebra structure,
i.e.~$A_\infty$-structures of the form $m = \partial + \mu + m_3 +
m_4 + \cdots$. Let us write $a = m_3 + m_4 + \cdots$. 

\medskip
Recall that $m = \partial + \mu + a$ is an $A_\infty$-structure 
if and only if $a$ satisfies the Maurer-Cartan formula and that 
such $a$ are called \emph{twisting cochains}. In this classical 
case the Maurer-Cartan formula reads
\[ -D(a)=\frac{1}{2}[a,a] \] if 2 is invertible in $k$ or,
equivalently, $-D(a)=a \circ a$ where $\circ$ denotes the composition 
product, see e.g. \cite[Section 2]{FiaPen02} and (\ref{derivedmaurercartan}). 

\begin{lem}\label{perturb} Let $A$ be a dga with differential $\partial$ and multiplication $\mu$,
and let $a$ be a twisting cochain. Further, for $n\ge 3$, let either $p \in
C^{n,1-n}(A,A)$ with $d^h(p)=[\mu,p]=0$ or  $p \in C^{n-1,2-n}(A,A)$
with $d^v(p)=[\partial,p]=0$. Then there is a twisting
cochain $\overline{a}$ such that
\begin{itemize}
\item the $A_\infty$-structures $\partial+\mu+\overline{a}$ and $\partial+\mu+a$
are quasi-isomorphic,
\item $\overline{a}_i = a_i$ for $i \le n-1$,
\item $\overline{a}_n= a_n - D(p)$. \qed
\end{itemize}
\end{lem}

We omit the proof since it is very similar to that of Lemma~\ref{derivedperturb}.
For the case where $A$ is a graded algebra rather than a dga, the
analogous result is mentioned without proof in~\cite[Section 4]{Kad88}.

With the help of Lemma~\ref{perturb}, we can now prove the
sufficient condition for a unique $A_\infty$-structure on a dga $A$
extending the existing differential and multiplication. 
This is only a minor generalization of Kadeishvili's classical result~\cite[Theorem 1]{Kad88} in the zero differential case, but we have not
been able to find a reference.

\medskip
To formulate the uniqueness results of this section and Section~\ref{sec:massey} 
we have to look deeper into the grading of the Hochschild cohomology of $A_\infty$-algebras 
and the internal grading of representing cocycles. 
An element of
$\hh^n(A,A)$ can be non-uniquely expressed as
\[
[x]=[x_0 + x_1 + x_2 + \cdots] \quad\mbox{with}\quad x_i \in
C^{i,n-i}(A,A).
\]
However, while the sum of the $x_i$ is a cocycle the individual
summands are not necessarily cocycles themselves. So generally we do
not get a decomposition of $\hh^n(A,A)$ as $\prod_i
\hh^{i,n-i}(A,A)$. To keep track of the internal degrees we
introduce a decreasing filtration on $\hh^*(A,A)$.

\begin{ddd}\label{def:filtration} For an $A_\infty$-algebra $A$, let
\[
F^k \hh^n(A,A) = \{ [x] \in \hh^n(A,A) \,|\, x \in
\prod\limits_{i \ge k} C^{i,n-i}(A,A) \}.
\]
\end{ddd}
This means that $F^k \hh^n(A,A)$ consists of all those elements of
$\hh^n(A,A)$ whose representing cocycles can be written as a sum of $x_i
\in C^{i,n-i}(A,A)$ with $i \ge k$.

Note that in the case of a bidga the filtration $F^*$ given in Definition
\ref{def:filtration} agrees with the usual filtration arising from
the column-wise filtration on the bicomplex, see e.g.~\cite[2.2 and
2.4]{Mcc01}.

\begin{theorem}
\label{th:HHcriterion}
Let $A$ be a dga with differential $\partial$ and multiplication $\mu$. If
\[
F^3 \hh^{2}(A,A)=0,
\]
then any $A_\infty$-structure on $A$ with
$m_1=\partial$ and $m_2 = \mu$ is quasi-isomorphic to $\partial+\mu$.
\end{theorem}

\begin{proof}
Let $a$ be a twisting cochain. Assuming that there is a $k \ge 3$
such that $a_i = 0$ for $i < k$, we are going to show that there is
a twisting cochain $\overline{a}$ that is equivalent to $a$ and
satisfies $\overline{a}_i=0$ for $i \le k$, i.e.~we are killing off
the bottom summand.  By induction, it follows that $a$ is equivalent
to zero.

So let $a$ now be a twisting cochain such that there is a $k \ge 3$
with $a_i = 0$ for $i < k$. Considering the Maurer-Cartan equation
\[
-D(a)=a \circ a
\]
in bidegrees $(k+1,2-k)$ and $(k,3-k)$, we see that $D(a_k)=0$ for
degree reasons, so $a_k$ is a cocycle and $[a_k] \in F^k
\hh^2(A,A)$. Since $F^k \hh^{2}(A,A)=0$, $a_k$ also has to be a
coboundary, i.e.~there is a cochain $p$ in total degree 1 with
$D(p)= a_k$. This $p$ is the sum of two cochains $p_1$ and $p_2$
with $p_1 \in C^{k,1-k}(A,A)$ and  $p_2 \in C^{k-1,2-k}(A,A)$. We
have $d^v(p_1) + d^h(p_2) = a_k$ and $d^h(p_1) = d^v(p_2) = 0$ for
degree reasons.

\medskip
$\xymatrix{
C^{k-1,3-k}(A,A) & C^{k,3-k}(A,A) & \cdots \\
C^{k-1,2-k}(A,A) \ar[u]^{d^v} \ar[r]^{d^h}  & C^{k,2-k}(A,A) \ar[u]^{d^v}  \ar[r]^{d^h} & C^{k+1,2-k}(A,A) \\
\cdots & C^{k,1-k}(A,A) \ar[u]^{d^v}  \ar[r]^{d^h}  & C^{k+1,1-k}(A,A)   \\
}$

\medskip
Applying Lemma \ref{perturb} for $p_1$ and $p_2$, we obtain that there is a twisting cochain $\overline{a}$ quasi-isomorphic to $a$ with $\overline{a}_i = 0 $ for $i < k$ and $\overline{a}_k=a_k - D(p) = 0$, which completes our proof.
\end{proof}

\begin{ex}
Consider the dga over the $p$-local integers
\[
A = \mathbb{Z}_{(p)}[x] \otimes \Lambda_{\mathbb{Z}_{(p)}}(e) \big/
(x^m, x^{m-1}e), \quad\quad \partial(x)=pe, \quad\quad |e|=-(2p-3),
|x|=-(2p-2)
\]
\smallskip
where $m \ge 2$. We can compute its Hochschild cohomology as a dga
by applying the spectral sequence for the homology of the total
complex of a bicomplex \cite[2.15]{Mcc01}. Its $E_1$-term is the
Hochschild cohomology of $A$ as a graded algebra.

\medskip To obtain this, we note that for an $A$-bimodule $M$
\begin{multline}
\hh^{*,*}_{alg}(A,M) \cong \hh^{*,*}_{alg}(\mathbb{Z}_{(p)}[x]\big/
(x^m)\otimes \Lambda_{\mathbb{Z}_{(p)}}(e),M) \nonumber \\\cong
\hh^{*,*}_{alg}(\mathbb{Z}_{(p)}[x]\big/ (x^m),M) \otimes
\hh^{*,*}_{alg}( \Lambda_{\mathbb{Z}_{(p)}}(e),\mathbb{Z}_{(p)}).
\nonumber
\end{multline}
(Use \cite[XI.1]{CarEil56} for the second isomorphism. The first
follows from a change-of-rings spectral sequence, see \cite{Mcc01}.)
Computing each factor separately, we obtain
\[
\hh^{*,*}_{alg}(A,A)= \mathbb{Z}_{(p)}[f,\tau] \otimes
\Lambda_{\mathbb{Z}_{(p)}}(\sigma) \otimes A
\]
with $|f|=(1,-|e|)$, $|\tau|=(2,-m|x|)$ and $|\sigma|=(1,-|x|)$ for
$A$ viewed as a graded algebra.

\medskip Already at this $E_1$-stage we can read off that
$\hh^{n,2-n}_{alg}(A,A)=0$ for $n \ge 3$, so $F^3 \hh^{2}(A,A)=0$
for $A$ as a dga. Hence $\mu+\partial$ is the only
$A_\infty$-structure on $A$ with $m_1=\partial$ and $m_2=\mu$.

\medskip Also note that the homology of $A$ coincides with the stable
homotopy groups of the $K_{(p)}$-local sphere in a certain range,
i.e. $H^{-i}(A)=\pi_i(L_1 S^0)$ for $0\leq i \leq (m-1)(2p-2)-1$.

\end{ex}

Combining Kadeishvili's result on minimal models with Theorem~\ref{th:HHcriterion}, we recover the following
result which we already stated earlier as Corollary~\ref{uniquedga}.

\begin{ccc}\label{uniquedga2}
Let $A$ be a dga over a ground field and $H^*(A)$ its homology algebra. Suppose that
\[
\hh^{n,2-n}_{alg}(H^*(A),H^*(A))=0 \quad\quad\mbox{for } n \ge 3.
\]
Then $A$ is intrinsically formal.
\end{ccc}

\begin{proof}
We apply Theorem~\ref{th:HHcriterion} to $H^*(A)$ with the trivial
differential to see that any $A_\infty$-structure on this is
quasi-isomorphic to the trivial one. So in particular the minimal
model is quasi-isomorphic to the trivial structure. But the minimal
model is quasi-isomorphic to $A$, so $A$ is formal.

Now given a dga $A'$ with $H^*(A')\cong H^*(A)$, the same argument shows
that $A'$ is also formal and thus that $A'$ is quasi-isomorphic to $A$.
\end{proof}

We note that the corollary follows from the special case of
Theorem~\ref{th:HHcriterion} where the dga has trivial differential.

\section{Massey products}
\label{sec:massey}

Massey products provide some very useful additional structure when
studying differential graded algebras and their homology. They are
closely related to Toda brackets in triangulated categories which
have strong applications in homotopy theory. Here we explain the
relationship between Massey products and the $m_3$ part of
$A_\infty$-structures; see also~\cite[Lemma 5.14]{BenKraSch05}.

In this section, $k$ denotes a field of characteristic not 2.

Let $A$ be a differential graded algebra and $\alpha_1, \alpha_2,
\alpha_3$ elements in the homology $H^*(A)$ such that $\alpha_1
\alpha_2 = 0$ and  $\alpha_2 \alpha_3 = 0$. That means that for
chosen representing cocycles $a_i$ of $\alpha_i$ there is an element
$u_i$ such that $d(u_i)= (-1)^{1+|a_{i}|}a_i a_{i+1}$. With those
elements, one can now define the Massey product of $\alpha_1,
\alpha_2$ and $\alpha_3$ as follows.

\begin{ddd}\label{definitionmasseyproduct}
Let $\alpha_1$, $\alpha_2$ and $\alpha_3$ be as above. 
Then the {\it Massey product} 
$\left< \alpha_1,\alpha_2,\alpha_3 \right> \subset H^{|a_1|+|a_2|+|a_3|-1}(A)$ 
is defined as the set of homology classes of the elements
\[
(-1)^{1+|a_1|}a_1 u_2 + (-1)^{1+|u_1|} u_1 a_3
\]
ranging over all possible choices of representing cocycles $a_i$ 
of the $\alpha_i$ and $u_i$ such that $d(u_i)= (-1)^{1+|a_{i}|}a_i a_{i+1}$.
\end{ddd}

Note that the Massey product $\left< \alpha_1,\alpha_2,\alpha_3
\right>$ is a set rather than an element as the choices one makes
can be altered by appropriate cocycles. Hence, if one fixes any $x$
in the Massey product, for any other $x'$ in the Massey product
there is a $y \in \alpha_1 H^{|\alpha_3|+|\alpha_2|-1}(A) \oplus
H^{|\alpha_2|+|\alpha_1|-1}(A) \alpha_3$ such that $x' = x + y$. The
group
\[\alpha_1 H^{|\alpha_3|+|\alpha_2|-1}(A) \oplus  H^{|\alpha_2|+|\alpha_1|-1}(A)  \alpha_3 \]
is called the {\it indeterminacy} of $\left<
\alpha_1,\alpha_2,\alpha_3 \right>$. So a Massey product consists of
only one element if and only if its indeterminacy is zero. For more
details on Massey products, see e.g. \cite[A.1.4]{Rav86}.

\begin{ex} Let $k$ be a field of characteristic different from $2$. 
Consider the following noncommutative differential graded algebra
\[
A = k\left<x,y\right> \big/ (x^3, y^2, xy=-yx), \quad \quad
\partial(x)=0, \partial(y)=x^2, \quad \quad |x|=2, |y|=3.
\]
Its homology has a copy of $k$ in degrees $0, 2, 5$ and $7$ and zero
elsewhere. Let $[x]$ and $[xy]$ denote the homology classes of $x$ and
$xy$ respectively. Then
\[
2[xy] = \left< [x],[x],[x] \right> \in H^5(A),
\]
the indeterminacy being zero for degree reasons.
\end{ex}

\begin{ex}
The dga
\[
A = \mathbb{Z}_{(p)}[x] \otimes \Lambda_{\mathbb{Z}_{(p)}}(e) \big/
(x^m, x^{m-1}e), \quad\quad \partial(x)=pe, \quad\quad |e|=-(2p-3),
|x|=-(2p-2)
\]
considered in the previous section has nontrivial Massey products.
Take $a_k$ to be an order $p$ element in $H^{-(2p-2)k+1}(A)$. Then
\[
\left<a_i,p,a_j\right>=a_{i+j}.
\]
This is related to the Toda bracket relation $\left<\alpha_i,p,\alpha_j\right>$ 
in the homotopy groups of the $K_{(p)}$-local sphere $\pi_*L_1 S^0$.
\end{ex}

\medskip

In the context of $A_\infty$-algebras, Massey products can be
reformulated using minimal models which were introduced in the
previous section. We quote the following well-known result (see also
\cite[Lemma 5.14]{BenKraSch05}).

\begin{lem}\label{m3lemma}
Let $A$ be a dga and $H^*(A)$ its minimal model with
$A_\infty$-structure $m$. Let $\alpha_1,\alpha_2, \alpha_3\in
H^*(A)$. If the Massey product $\left< \alpha_1,\alpha_2,\alpha_3
\right>$ is defined in $H^*(A)$, then
$$(-1)^{|\alpha_1|+|\alpha_2|+1}m_3(\alpha_1\otimes\alpha_2\otimes\alpha_3)
\in \left< \alpha_1,\alpha_2,\alpha_3 \right>.$$\qed 
\end{lem}

Hence, if $A$ and $B$ are differential graded algebras with
isomorphic homology algebras $H^*(A)$ and $H^*(B)$, then they have
the same Massey products if the $A_\infty$-structures of the minimal
models have identical $m_3$. (The converse is not necessarily true,
see the discussion at the end of this section.)

\begin{theorem}\label{uniquenessmassey}
Let $A$ be a dga whose minimal model $H^*(A)$ satisfies $m_i = 0$
for $i \neq 2, 3$ and let $\overline{m}$ be an $A_\infty$-structure
on $H^*(A)$ with $\overline{m}_2 = m_2$ and $\overline{m}_3=m_3$. If
$F^4 \hh^{2}(H^*(A),H^*(A))=0$, then $\overline{m}$
and $m$ are quasi-isomorphic.
\end{theorem}

\begin{proof}
The proof is extremely similar to the proof of Theorem \ref{th:HHcriterion}.
The differential in the Hochschild complex for $H^*(A)$ is $D= D_2 + D_3$ with
\[
D_2 = [m_2,-]: C^{n,k}(H^*(A),H^*(A)) \longrightarrow
C^{n+1,k}(H^*(A),H^*(A))
\]
and
\[
D_3 = [m_3,-]: C^{n,k}(H^*(A),H^*(A)) \longrightarrow
C^{n+2,k-1}(H^*(A),H^*(A)).
\]
Assume there is an $A_\infty$-structure $\overline{m}$ on $H^*(A)$
with $\overline{m} = m_2 + m_3 + a_4 + a_5 + \cdots$. Let $a= a_4 +
a_5 + \cdots$. Because $m= m_2 + m_3$ is an $A_\infty$-structure on
the minimal model by assumption, we know that $a$ is a twisting
cochain, i.e.~$a$ satisfies the Maurer-Cartan equation. Again, for
degree reasons $D(a_4)= 0$ and so there is
 $p_2 \in C^{3,-2}(H^*(A),H^*(A))$ and
 $p_3 \in C^{2,-1}(H^*(A),H^*(A))$ with $D_2(p_2)+D_3(p_3)=a_4$ and $D_3(p_2)=D_2(p_3)=0$.
 The analogue of Proposition \ref{perturb} also holds in this case: for any
 $p \in C^{n,1-n}(H^*(A),H^*(A))$ with $D_3(p)=0$
 or $p\in C^{n+1, -n}(H^*(A), H^*(A))$ with $D_2(p)=0$,
 there is a twisting cochain $\overline{a}=\overline{a}_4 + \overline{a}_5 + \cdots$ such that
\begin{itemize}
\item $\overline{a}$ is equivalent to $a$,
\item $\overline{a}_k = {a}_k$ for $k\le n$,
\item $\overline{a}_{n+1} = a_{n+1} - D(p)$.
\end{itemize}
The rest of the proof follows the same steps as the proof of
Theorem \ref{th:HHcriterion}.
\end{proof}

Of course one would like to apply this theorem to a minimal model
$(H^*(A),m=m_2+m_3)$ of a dga $A$ to obtain a uniqueness result
analogous to Corollary \ref{uniquedga} and conclude that the
vanishing of the right Hochschild cohomology groups implies that $A$
is the only dga up to quasi-isomorphism with the given homology and Massey products.

This does not quite work- to give the same Massey products on
minimal models of dgas with the same homology algebras, $m_3$ only
needs to agree on triples $(a,b,c)$ with $ab=0=bc.$ For example, in
\cite[Example 5.15 and Proposition 5.16]{BenKraSch05} Benson, Krause
and Schwede constructed an example of a dga with trivial Massey
products but nontrivial $m_3$.

It would also be interesting to study the implication of Massey products regarding uniqueness criteria in the derived case.

\begin{appendix}

\section{Signs in the Lie bracket}

In this appendix we verify the signs appearing in the Lie
bracket of Section~\ref{trigraded}. The special case where
$k=l=0$ recovers the signs in Section~\ref{liealgebrastructure}.

\begin{lem}
In the context of Section \ref{trigraded},
\begin{align}
    [f,g] & :=  \sigma^{-1}\bc{\sigma(f),\sigma(g)} \nonumber\\
            & =\sum\limits_{v=0}^{n-1} (-1)^{(n-1)(m-1)+v(m-1)+j(n-1)} f(1^{\otimes
                v}\otimes g \otimes 1^{\otimes n-v-1}) \nonumber \\
            & -(-1)^{\left<f,g\right>} \sum\limits_{v=0}^{m-1}
                (-1)^{(m-1)(n-1)+v(n-1)+i(m-1)} g(1^{\otimes v}\otimes f \otimes
                1^{\otimes m-v-1}) \nonumber \\
        & \qquad\in C^{n+m-1,i+j}_{k+l}(A,A)
            \nonumber
\end{align}
for $f \in C^{n,i}_k(A,A)$ and $g \in C^{m,j}_l(A,A)$. Here,
$\left<f,g\right> := (n+i-1)(m+j-1)+kl$.
\end{lem}

\begin{proof}
Throughout this proof, by $\circ$, we mean the actual composition of morphisms rather than
the previously used composition product.

The signs arise from the Koszul sign rule for interchanging
morphisms. For morphisms $f, g, h$ and $u$, we have
\[
(f \otimes g) \circ (h \otimes u) = (-1)^{is+jt} (f \circ h) \otimes
(g \circ u)
\]
with $g$ having internal bidegree $(i,j)$ and $h$ having internal
bidegree $(s,t)$.

\medskip
We then obtain
\begin{align}
\sigma^{-1} \bc{\sigma(f),\sigma(g)} & = \sigma^{-1} \Big(
\sum\limits_{v=0}^{n-1} \sigma(f)(1^{\otimes v} \otimes \sigma(g)
\otimes 1^{\otimes n-v-1})  \Big) \nonumber \\ &
-(-1)^{\left<f,g\right>} \sigma^{-1} \Big( \sum\limits_{v=0}^{m-1}
\sigma(g)(1^{\otimes v} \otimes \sigma(f) \otimes 1^{\otimes m-v-1})
\Big). \nonumber
\end{align}
For reasons of symmetry and linearity we are only going to
explicitly compute
\[
\sigma^{-1}\Big(\sigma(f)(1^{\otimes v} \otimes \sigma(g) \otimes
1^{\otimes n-v-1})\Big).
\]
Up to sign, this is $f(1^{\otimes v} \otimes g\otimes 1^{\otimes
n-v-1})$ and we now calculate the sign.
\smallskip

The term $ \sigma(f)(1^{\otimes v} \otimes \sigma(g) \otimes
1^{\otimes n-v-1})$ lies in $C^{n+m-1,i+j+n+m-2}_{k+l}(S(A),S(A))$,
so
\begin{align*}
  &\sigma^{-1}\Big(\sigma(f)(1^{\otimes v} \otimes \sigma(g) \otimes 1^{\otimes n-v-1}) \Big)
  \nonumber \\
  &= (-1)^{i+j+n+m+k+l+\binom{n+m-1}{2}} S^{-1} \circ \Big(\sigma(f)(1^{\otimes v} \otimes
            \sigma(g) \otimes 1^{\otimes n-v-1})\Big) \circ S^{\otimes n+m-1}
            \nonumber \\
   &= (-1)^{\binom{n+m-1}{2}}S^{-1} \circ S \circ f \circ
            (S^{-1})^{\otimes n} \circ \Big(1^{\otimes v} \otimes (S \circ g
            \circ (S^{-1})^{\otimes m}) \otimes 1^{\otimes n-v-1}\Big) \circ
            S^{\otimes n+m-1} \nonumber \\
    &= (-1)^{\binom{n+m-1}{2}}  f \circ \Big((S^{-1})^{\otimes v} \otimes
            S^{-1} \otimes (S^{-1})^{\otimes n-v-1}\Big) \\
    &\qquad\qquad\qquad\qquad\qquad\qquad\qquad\qquad\circ \Big(1^{\otimes
            v} \otimes (S \circ g \circ (S^{-1})^{\otimes m}) \otimes 1^{\otimes
            n-v-1}\Big) \circ S^{\otimes n+m-1}. \nonumber
\end{align*}
In the next step we are
obtaining a new sign $(-1)^{(n-v-1)(j+m-1)}$ by interchanging
$(S^{-1})^{\otimes n-v-1}$ with $1^{\otimes v} \otimes (S \circ g
\circ (S^{-1})^{\otimes m})$. Interchanging $S^{-1}$ and $1^{\otimes
v}$ does not introduce any new signs, so we continue with
\begin{multline}
(-1)^{\binom{n+m-1}{2} + (n-v-1)(j+m-1)}  f \circ
\Big((S^{-1})^{\otimes v} \otimes ( S^{-1} \circ S \circ g \circ
(S^{-1})^{\otimes m}) \otimes (S^{-1})^{\otimes n-v-1}\Big) \circ
S^{\otimes n+m-1} \nonumber
\\
= (-1)^{\binom{n+m-1}{2} + (n-v-1)(j+m-1)} \nonumber \\  f \circ
\Big((S^{-1})^{\otimes v} \otimes  g \circ (S^{-1})^{\otimes m}
\otimes (S^{-1})^{\otimes n-v-1}\Big) \circ \Big(S^{\otimes
v}\otimes S^{\otimes m} \otimes S^{\otimes n-v-1} \Big).\nonumber
\end{multline}
Since $(S^{-1})^{\otimes a} = (-1)^{\binom{a}{2}} (S^{\otimes
a})^{-1}$, we continue with
\begin{multline}
(-1)^{\binom{n+m-1}{2} + \binom{n-v-1}{2} + \binom{m}{2} +
\binom{v}{2} + (n-v-1)(j+m-1)} \nonumber \\ f \circ \Big((S^{\otimes
v})^{-1} \otimes  g \circ (S^{\otimes m})^{-1} \otimes (S^{\otimes
n-v-1})^{-1}\Big) \circ \Big(S^{\otimes v}\otimes S^{\otimes m}
\otimes S^{\otimes n-v-1} \Big).\nonumber
\end{multline}
We have that
\[
\binom{n+m-1}{2} + \binom{n-v-1}{2} + \binom{m}{2} + \binom{v}{2}
\equiv (n-1)(v+m) + v \pmod 2
\]
so we can simplify the sign in the above expression to give
\begin{multline}
(-1)^{ (n-1)(v+m) + v + (n-v-1)(j+m-1)} \nonumber \\ f \circ
\Big((S^{\otimes v})^{-1} \otimes  g \circ (S^{\otimes m})^{-1}
\otimes (S^{\otimes n-v-1})^{-1}\Big) \circ \Big(S^{\otimes
v}\otimes S^{\otimes m} \otimes S^{\otimes n-v-1} \Big).\nonumber
\end{multline}
We then interchange $S^{\otimes v}$ with $g \circ (S^{\otimes
m})^{-1} \otimes (S^{\otimes n-v-1})^{-1}$ which in addition gives
us the new sign $(-1)^{(j+m+n-v-1)v}$, so we have
\begin{multline}
(-1)^{ (n-1)(v+m) + v + (n-v-1)(j+m-1)+(j+m+n-v-1)v} \nonumber \\ f
\circ \Big(1^{\otimes v} \otimes  g \circ (S^{\otimes m})^{-1}
\otimes (S^{\otimes n-v-1})^{-1}\Big) \circ \Big(1^{\otimes
v}\otimes S^{\otimes m} \otimes S^{\otimes n-v-1} \Big).\nonumber
\end{multline}
Finally, we add to the sign by interchanging $S^{\otimes m}$ with
$(S^{\otimes n-v-1})^{-1}$, so we end up with
\begin{multline}
(-1)^{ (n-1)(v+m) + v + (n-v-1)(j+m-1)+(j+m+n-v-1)v+ m(n-v-1)} f
\circ (1^{\otimes v} \otimes  g  \otimes 1^{\otimes n-v-1}).
\nonumber
\end{multline}
We can then simplify the above sign to
\[
(-1)^{(n-1)(m-1)+v(m-1)+(n-1)j}
\]
which proves our claim.
\end{proof}

\end{appendix}

\end{document}